# Queues with Censored Demand and Autoregressive Net Input

Kerry Fendick
fendick@att.net

*Abstract*

We develop methods for simulating a queue with an autoregressive net-input process and for recovering characteristics of such a net-input process from samples of queue lengths. We apply these methods to the problem of estimating the censored (unsatisfied) demand for the queue's content and show how to model a queue for which the censoring of demand is graduated in a neighborhood of the queue's zero lower bound. As an example, we estimate the monthly unsatisfied demand for U.S. nonfarm jobs based on samples of job openings through a period including the last two recessions.

## 1  Introduction

This paper is a sequel to Fendick [1], which derived properties of Gaussian Markov processes with stationary increments (GMSI processes) and of queues driven by them. GMSI processes can exhibit positively correlated, negatively correlated, or uncorrelated increments. In the present paper, we introduce methods for estimating parameters and conditional expectations of these and related models using time-series data.

Here, we model a continuously evolving queue that is observed only periodically. In the case on which we focus most of our attention, the queue lengths are observed, the queue's net-input increments are not observed, but there may exist samples of other observed processes predictive of the net-input increments. Our model is characterized by two distributions: one characterizing the autoregressive dependencies of the net-input process without regard to the queue's dynamics; the other characterizing the queue's boundary behavior without regard to its autoregressive dependencies. We derive the second of those distributions originally assuming that demand for the queue's content is censored only when the queue is empty but we later allow the censoring of demand also to occur at positive queue lengths in a neighborhood of zero. We apply the model to estimate the expected amount of censored demand in each period.

In characterizing economic entities, the field of accounting distinguishes between flow and stock variables. Flow variables represent a change in a quantity over an interval, whereas stock variables represent accumulated quantities at snapshots in time. The results of this paper are useful for modeling stock variables as queue lengths and estimating the censored demand for the quantity represented by the stock variable.

In many examples of practical interest, stock variables are defined to be non-negative and can be modeled as the length of a queue that increases according to a cumulative supply process and decreases according to a cumulative demand process. In such examples, it is often not easy to discern the total demand for the quantity represented by the stock variable (that is, by the amount of stock). In Section 6, for example, we will model U.S. non-farm job openings - a stock variable - as a sequence of queue lengths. Although monthly data is available on the number of unemployed seeking work, it is a matter of debate among economists how much unsatisfied demand for jobs in any given month is constrained by the level of job openings. A portion of the unemployed may no longer be seeking job energetically. For those who are, the available openings may not be suitable: for some, the available openings may be unappealing relative to their past jobs; for others, the available openings may require new skills that have not yet been acquired. In such cases there is unsatisfied demand for jobs, but not necessarily for the openings that are available; and more of the same openings for a given month would not necessarily result in more hiring that month.

We will use the term "censored demand" to refer to unsatisfied demand attributable to insufficient stock. Another example of a stock variable that can be modeled as a length of a queue is the inventory of a firm (or more precisely, the cost of the goods in a firm's inventory). When inventory levels for a firm decrease, lead times for delivery of purchased items may increase and demand for the firm's goods may be censored from customers who take their business elsewhere as a result. The censored demand each period then translates to additional sales that would have occurred had the firm maintained more of the same inventory, so that estimates of it would be useful in evaluating economic tradeoffs of holding more or less inventory.

We model a stock variable as a queue under the assumption that the unobserved net-input to the queue (the difference between the queue's cumulative supply and cumulative total demand) has stationary increments that are Gaussian and Markov. Gaussian approximations are sometimes criticized for their inability to account for non-negativity constraints, but we exploit that limitation here. When queue lengths are long enough, their change over an interval equals the change in the net-input process. Differences between the observed characteristics of the queue-length process and the assumed characteristic of the net-input process reveal the expected magnitude of censored demand in each period.

## 1.1 Organization of Paper

The remainder of the paper is organized as follows. Section 1.2 provides definitions for the model we will analyze, while leaving most probabilistic assumptions for later. Section 1.3 compares and contrasts the work here with related literature. Section 2 reviews the definition and properties of GMSI processes and extends them by deriving an alternative parameterization that facilitates parameter estimation. The properties of GMSI processes provide the basis for a queuing model that we develop in Section 3. Section 3.1 contains an algorithm defining that model, and Section 3.2 presents simulation results showing how changes to the model's parameters affect the lengths of queues and busy periods.

Section 4 shows how to estimate the parameters of the queuing model given samples from a queue's net-input process. These methods are used as building blocks for what comes later, but are potentially relevant in themselves for studying communications and service systems, since measurements of the net-input process are often available for those applications. Section 4.1 provides details on consistent estimators of covariances and autogressive coefficients, and Section 4.2 tests the methods using samples from simulations. Section 5 then shows how to estimate model parameters and conditional expectations when samples of queue lengths are available, but samples of the net-input process are not. Using samples from simulations defined in Section 3, Section 5.1 tests the methods for estimating parameters, and Section 5.2 examines properties of estimated conditional expectations.

Section 6 generalizes the models for queues at which demand is censored not only at the queue's zero lower bound but also at positive queue lengths in some neighborhood of zero. Section 6.2 presents examples of simulations of such queues, and Section 6.3 presents an example of estimating parameters of such queues based on queue lengths.

Section 7 is an extended example applying the results of this paper to jobs data as discussed earlier. There, we discuss techniques for transforming the raw jobs data so that it conforms more closely to model assumptions. Similar considerations would apply to other examples in which stock variables are modeled as queues.

### 1.2   A Continuous-time Queue at Discrete Epochs

Although a stock variable, by definition, describes a quantity at a snapshot in time, the quantity itself usually changes more or less continuously between snapshots. We will model the evolution of such a quantity as a continuous-time queue $\{Q(t): t \geq 0\}$ satisfying

$$Q(t) = Q(0) + Z(t) + L(t) \geq 0 \; for \; t \geq 0 \quad (1.1)$$

where

$L(\cdot)$ is a non-decreasing continuous process with $L(0) = 0$ \quad (1.2)

and

$$L(\cdot) \text{ increases only when } Q(\cdot) = 0. \qquad (1.3)$$

At time $t$, the random variable $Q(t)$ has an interpretation as a queue's content, and the random variable $Z(t)$ as the queue's cumulative net-input, that is, the cumulative supply of content to the queue minus the cumulative demand for content from the queue. The random variable $L(t)$ then has the interpretation as the cumulative censored demand, that is, the cumulative demand for the queue's content that is unsatisfied over the interval $[0,t]$. The functional defined by (1.1)-(1.3) mapping $Z$ to $Q$ is variously known as the *one-dimensional reflection map*, the *one-sided regulator*, and the *Skorokhod map*; for background see Section 5.2.2 and Chapter 14 of Whitt [2] and Section 2, Chapter 2 of Harrison [3].

Next let
$$q_n \equiv Q(n) \ for \ n = 0, \dots, N$$
$$z_{n,1} \equiv Z(n) - Z(n-1) \ for \ n = 1, \dots, N \qquad (1.4)$$

and
$$l_n \equiv L(n) - L(n-1) \ for \ n = 1, \dots, N$$

Then (1.1) and (1.2) imply that
$$q_n = q_{n-1} + z_{n,1} + l_n \geq 0 \ for \ n = 1, \dots, N \qquad (1.5)$$

where
$$l_n \geq 0 \text{ for } n = 1, \dots, N. \qquad (1.6)$$

For period $n$, the random variable $z_{n,1}$ is the net-input increment, and $l_n$ the censored demand. In accounting terms, the $q_n$'s are stock variables, and the $z_{n,1}$'s and $l_n$'s are flow variables. In the general case, we will also assume that there exists other flow variables $z_{n,i}$ for $i = 2, \dots, K$ and $n = 1, \dots, N$ that are predictive of the net-input increments. (In the case for which $K = 1$, these other increments are not present in the model.)

One goal of this paper is to estimate the expected value of the unobserved censored demand $l_n$ for each period $n$ conditional on the observed queue lengths. The estimates will depend implicitly on the continuous-time dynamics of the underlying processes in (1.1)-(1.3). In solving that problem, we assume that samples of net-input increments $\{z_{n,1}\}$ are unavailable since the total demand is a constituent of the net-input process; and the censored demand is a constituent of the total demand.

The assumption in (1.3) is an idealization for the quantities we are interested in modeling. Examples of stock variables in economics and finance often have operational lower bound that are non-zero. The characteristic of having a non-zero operational lower bound distinguishes stock variables from other quantities more commonly modeled as queue lengths. As examples, the number of customers waiting for a typical bank teller is equal to zero part of day, but the reported number of job openings for the US economy is never equal zero, even during the deepest recessions.

In Section 7, we provide an example in which we transform data to conform to (1.3). Even after such a transformation, (1.3) is often at best a rough approximation to how demand for the types of quantities represented by a stock variable is often censored. As another example, if the gasoline inventory of an oil and gas company were to become low, shortages might occur first in some regions, so that censorship of demand would occur in a graduated fashion while the total inventory was still positive. Modeling the inventory of gasoline inventories of different regions separately might improve the accuracy of (1.3) as an approximation, but shortages might still occur first at some gas stations; so that censorship of demand would still be graduated within a region.

We will say that a queue has a sharp boundary if it satisfies (1.1)-(1.3) and a graduated boundary if it satisfies (1.1) and (1.2) but not (1.3). In Section 6, we show how to model a queue with a graduated boundary. It is well known that the queue-length process $Q$ in (1.1) is uniquely specified when the net-input process $Z$ is continuous and the boundary is sharp; see for example Proposition 3 on page 19 of Harrison [3]. The queue-length process is not uniquely specified when the boundary is graduated as defined to this point in the paper. Nevertheless, we show in Section 6 how additional assumptions relating the net-input process to GMSI processes lead to a natural approach for specifying models with graduated boundaries.

## 1.3 Related Work

Queuing models have been widely applied in communications, computer design, traffic analysis, manufacturing, service industries, and finance. Kleinrock [4], Newell [5], and Whitt [2] provide numerous examples. Queuing models have been widely applied to model inventories; see for example Harrison [3] for development of inventory models for which the net-input processes is modeled as Brownian motion with drift, a special case of the model developed here. Queuing models have been recently applied by Cont, Stoikov, & Talreja [6] for order-book dynamics of financial exchanges.

The problem of censored demand has been studied in the context of inventory control; see Section 1.2 of Huh and Rusmevichientong [7] for a recent literature survey and Conrad [8] for original work in this area. Previous work on the problem of censored demand for inventory has assumed that the demand in different periods is independent and identically distributed and that the realized demand (in the form of sales) is observed. In applying the results here to an inventory, only the inventory levels themselves need be observed, and the unobserved demand in different periods can be autocorrelated (as reflected by autocorrelations in the net-input process). Most previous models for estimating censored demand have addressed aspects of inventory control, such as stocking strategies and perishability, not addressed here.

Section 5 here provides an example of estimating parameters of a stochastic process based on a realization that is only partially observed. It is an example of a problem addressed in Schnedler [9] in much greater generality. Schnedler's solution through likelihood estimation assumes that the joint distribution of all variables comprising the realization are known (though its parameters are not). In the context of Section 1.2, that would require knowing the joint distribution of

$$[q_0, \ldots q_N, z_{1,1}, \ldots, z_{N,1}, \ldots, z_{1,K}, \ldots, z_{N,K}]$$

except for some parameters. Instead, we estimate the parameters assuming that certain conditional distributions involving these random variables are known. We make the further simplifying assumption that the only important dependences on past samples are over a finite history. Different assumptions for the length of that history result in different models, and we select the model achieving the best fit of observed data.

The Expectation Maximization (EM) Algorithm, as described in Section 10.4.6 on page 375 of Amemiya [10], is a general approach to estimating the parameters of models with latent (unobserved) variables. The algorithm that we develop in Section 5 is similar to the EM algorithm in that both can be described as iteratively generating samples of the unobserved variables based on current estimates of the model's parameters and then updating those estimates based on those samples. For the EM algorithm, the samples of the unobserved variables are taken to be their conditional expectations as a function of the current parameter estimates. For the algorithm in Section 5 here, those samples are instead drawn from conditional distributions that depend on the current parameter estimates. We call the algorithm in Section 5 the Distribution Minimization (DM) algorithm to highlight its similarities and differences relative to the EM algorithm. The parameter estimates obtained from the EM algorithm are the estimates to which the EM algorithm converges; the parameter estimates obtained from the DM algorithm are the averages of the corresponding estimates obtained in the algorithm's steps. The EM algorithm incorporates maximum likelihood estimation; the DM algorithm solves a least-squares minimization.

Stochastic processes can exhibit different degrees of variability when measured over different time scales. Indices of dispersion are a means of quantifying covariances as a function of time scale; see Cox & Lewis [11] and Whitt [2] for background. Central limit theorems, as covered by Whitt [2], can often describe complex stochastic processes through a small number of parameters because the relevant time scales for the processes become large in those limits. Estimating performance measures for a system not operating in such a limiting regime explicitly or implicitly requires determining its relevant history (or time scale). Here, we map the problem of determining a relevant history to the model selection problem of regression analysis. The covariances we use in fitting the models convey similar information to the values of indices of dispersion at a given time scale, as both describe partial sums of increments. But, whereas that information was used by past studies (as surveyed in Section 9.6 of Whitt [2] ) to fit models with independent increments, we use it to fit autoregressive models with autocorrelated increments.

Queuing systems for which (1.1) and (1.2) hold but not (1.3) have been studied previously assuming explicit policies controlling the process $L(\cdot)$; see for example Chapter 6 of Harrison [3] and references cited there. In our model of systems with graduated boundaries developed in Section 6, the behavior of $L(\cdot)$ is implicitly defined through a distribution function that can be used to tune the amount of demand censored away from the boundary at zero to fit data.

## 2  GMSI Processes

We begin by reviewing and extending results from Fendick [1] characterizing Gaussian Markov process with stationary increments. Let $\{X^{(i)}(t): 0 \leq t < \Delta\}$ for $i = 1, \ldots, K$ denote a real, continuous zero-mean Gaussian process with $X^{(i)}(0) = 0$ and covariance function

$$r^{(i)}(s,t) \equiv E[X^{(i)}(s)X^{(i)}(t)] = s(\theta^{(i)} - \tau^{(i)}t) \text{ on } 0 \leq s \leq t < \Delta \leq \infty, \qquad (2.1)$$

where $\theta^{(i)} > 0$ and $\tau^{(i)} \leq \theta^{(i)}/\Delta$ are constants. By Theorem 1 of Fendick [1], $X^{(i)}$ is a canonical representation of a zero-mean Gaussian Markov process with stationary increments and smooth covariance function. We call $X^{(i)}$ a $(\theta^{(i)}, \tau^{(i)})$ *GMSI process* on $[0, \Delta)$. By Corollary 2 of Fendick [1], a $(\theta^{(i)}, \tau^{(i)})$ *GMSI* process has negative autocorrelations when $\tau^{(i)}$ is positive, and positive autocorrelations when $\tau^{(i)}$ is negative.

We will further assume that the $X^{(i)}$'s are independent of one another, and let

$$X(t) \equiv \sum_{i=1}^{K} k_i X^{(i)}(t) \text{ on } 0 \leq t < \Delta \qquad (2.2)$$

where the $k_i$'s are constants. By Proposition 1 of Fendick, $X$ is itself a $(\theta, \tau)$ *GMSI process on* $[0, \Delta)$ with

$$\theta \equiv \sum_{i=1}^{K} k_i^2 \theta^{(i)} > 0 \text{ and } \tau \equiv \sum_{i=1}^{K} k_i^2 \tau^{(i)} \leq \frac{\theta}{\Delta}. \qquad (2.3)$$

so that

$$r(s,t) \equiv E[X(s)X(t)] = s(\theta - \tau t) \text{ on } 0 \leq s \leq t < \Delta \leq \infty. \qquad (2.4)$$

By Proposition 2 of Fendick [1],

$$P(X(M+1) \leq x | \ X(M) = \sum_{i=1}^{K} k_i x_i, \ X^{(i)}(M) = x_i \ for \ i = 2, \dots, K)$$

$$= P(X(M+1) \leq x | X^{(i)}(M) = x_i \ for \ i = 1, \dots, K)$$

$$= N\left(x; \sum_{i=1}^{K} \left(k_i x_i - \frac{k_i x_i \tau^{(i)}}{\theta^{(i)} - M\tau^{(i)}}\right), \theta - \sum_{i=1}^{K} k_i^2 \frac{\theta^{(i)} \tau^{(i)}}{\theta^{(i)} - M\tau^{(i)}}\right) \quad (2.5)$$

when $\Delta > M + 1$ where $dN(x; \mu, \sigma^2)/dx = exp(-(x-\mu)^2/(2\sigma^2))/\sqrt{2\pi\sigma^2}$ is the density of a normally distributed random variable with mean $\mu$ and variance $\sigma^2$.

We now express the $3K$ parameters of $X(\cdot)$ defined in (2.2), namely $\{k_i\}_{i=1}^{K}, \{\theta^{(i)}\}_{i=1}^{K}$, and $\{\tau^{(i)}\}_{i=1}^{K}$, as functions of other parameters that are easier to estimate. The vector $\left(X(M+1), X(M), X^{(2)}(M), \dots, X^{(K)}(M)\right)$ has a multivariate Gaussian distribution with a zero mean vector and a covariance matrix that we will denote by $\mathbf{\Sigma}$. The elements of $\mathbf{\Sigma}$ are functions of the aforementioned $3K$ parameters. If we partition $\mathbf{\Sigma}$ as

$$\mathbf{\Sigma} = \begin{bmatrix} \mathbf{A} & \mathbf{B} \\ \mathbf{B}^T & \mathbf{\Gamma} \end{bmatrix} \quad (2.6)$$

with respective sizes

$$\begin{bmatrix} 1 \times 1 & 1 \times K \\ K \times 1 & K \times K \end{bmatrix}$$

then it follows from well-known properties of the multivariate normal distribution that

$$E\left(X(M+1) | \ X(M) = \sum_{i=1}^{K} k_i x_i, \ X^{(i)}(M) = x_i \ for \ i = 2.,,. K\right)$$

$$= \mathbf{B}\mathbf{\Gamma}^{-1} \begin{bmatrix} \sum_{i=1}^{K} k_i x_i \\ x_2 \\ \vdots \\ x_K \end{bmatrix} \quad (2.7)$$

see for example P3 on page 197 of Goldberg [12].

The covariance matrix $\mathbf{\Sigma}$ depends on the parameters $k_1$, $\theta^{(1)}$, and $\tau^{(1)}$ only through the products $k_1 \times \theta^{(1)}$ and $k_1 \times \tau^{(1)}$, so we will let

$$k_1 = 1 \quad (2.8)$$

without loss of generality.

The conditional expectation in (2.7) is a linear combination of $\sum_{i=1}^{K} x_i, x_2, \dots, x_K$ with respective *autoregression coefficients* that we will denote by

$\xi_1, \xi_2, \ldots, \xi_K$. We obtain $K$ equations in the desired parameters through the identities

$$\langle \mathbf{B}\mathbf{\Gamma}^{-1}\rangle_i = \xi_i \text{ for } i = 1, \ldots, K \tag{2.9}$$

where $\langle \mathbf{y}\rangle_i$ denotes the $i^{th}$ component of the vector $\mathbf{y}$.

The matrix $\mathbf{\Gamma}$ in (2.6) is symmetric with upper triangular elements

$$\mathbf{\Gamma} = \begin{bmatrix} r(M,M) & k_2 r^{(2)}(M,M) & k_3 r^{(3)}(M,M) & \cdots & k_{K-1} r^{(K-1)}(M,M) & k_K r^{(K)}(M,M) \\ & r^{(2)}(M,M) & 0 & \cdots & 0 & 0 \\ & & r^{(3)}(M,M) & \ddots & 0 & 0 \\ & & & \ddots & 0 & 0 \\ & & & & r^{(K-1)}(M,M) & 0 \\ & & & & & r^{(K)}(M,M) \end{bmatrix} \tag{2.10}$$

where $r(\cdot,\cdot)$ is defined in (2.4), and $r^{(i)}(\cdot,\cdot)$ as in (2.1). Defining

$$\gamma_{1,1} \stackrel{\text{def}}{=} Var\left(X(M)\right), \quad \gamma_{i,i} \stackrel{\text{def}}{=} Var\left(X^{(i)}(M)\right) \text{ for } i = 2, \ldots, K \tag{2.11}$$

and

$$\gamma_{1,i} \stackrel{\text{def}}{=} Cov\left(X(M), X^{(i)}(M)\right) \text{ for } i = 2, \ldots, K. \tag{2.12}$$

we obtain the remaining $2K - 1$ equations in the desired parameters through the following identities corresponding to the non-zero elements of $\mathbf{\Gamma}$:

$$r(M,M) = \gamma_{1,1} \tag{2.13}$$

$$r^{(i)}(M,M) = \gamma_{i,i} \text{ for } i = 2, \ldots, K \text{ and} \tag{2.14}$$

and

$$k_i r^{(i)}(M,M) = \gamma_{1,i} \text{ for } i = 2, \ldots, K. \tag{2.15}$$

In Section 4, we will describe how to estimate the autoregressive coefficients $\xi_i$ and covariances $\gamma_{ij}$ from samples.

*Proposition 1: For given $K \geq 1$ and $M \geq 1$, if*

$$\xi_1 < 1 + \frac{1}{M} \tag{2.16}$$

*and*

$$\xi_1 + \xi_i \frac{\gamma_{i,i}}{\gamma_{1,i}} < 1 + \frac{1}{M} \text{ for } i = 2, \ldots, K, \tag{2.17}$$

*where*

$$\gamma_{i,i} > 0 \text{ for } i = 1, \ldots, K, \tag{2.18}$$

$$\gamma_{1,i} = \gamma_{i,1} \neq 0 \text{ for } i = 2, \ldots, K, \tag{2.19}$$

$$\gamma_{j,i} = 0 \text{ otherwise}, \tag{2.20}$$

and

$$\sum_{i=2}^{K} \frac{\gamma_{1,i}^2}{\gamma_{1,1}\gamma_{i,i}} \leq 1 \tag{2.21}$$

*then the equations (2.8), (2.9), and (2.13)-(2.15) are solved by*

$$k_i = \begin{cases} \gamma_{1,i}/\gamma_{i,i}, & i = 2, \ldots, K \\ 1, & i = 1 \end{cases} \tag{2.22}$$

$$\theta^{(i)} = \begin{cases} \gamma_{i,i}\left(1 + \frac{1}{M} - \left(\xi_1 + \xi_i \frac{\gamma_{i,i}}{\gamma_{1,i}}\right)\right), & i = 2, \ldots, K \\ \left(\gamma_{1,1} - \sum_{i=2}^{K} \gamma_{1,i}^2/\gamma_{i,i}\right)\left(1 + \frac{1}{M} - \xi_1\right), & i = 1 \end{cases} \tag{2.23}$$

$$\tau^{(i)} = \begin{cases} \frac{\gamma_{i,i}}{M}\left(1 - \left(\xi_1 + \xi_i \frac{\gamma_{i,i}}{\gamma_{1,i}}\right)\right), & i = 2, \ldots, K \\ \left(\gamma_{1,1} - \sum_{i=2}^{K} \gamma_{1,i}^2/\gamma_{i,i}\right)\frac{1 - \xi_1}{M}, & i = 1 \end{cases} \tag{2.24}$$

*where*

$$\theta^{(i)} > 0 \text{ and } \tau^{(i)} \leq \frac{\theta^{(i)}}{M+1} \text{ for } i = 1, \ldots, K. \tag{2.25}$$

*When $K = 1$, (2.22)-(2.24) reduce to $k_1 = 1$, $\theta^{(1)} = \gamma_{1,1}(1 + 1/M - \xi_1)$ and $\tau^{(1)} = \gamma_{1,1}(1 - \xi_1)/M$.*

Proof of Proposition 1:

The inequalities in (2.25) are easily verified. When (2.22)-(2.24) hold, it is also easily verified that (2.8) and (2.13)-(2.15) are satisfied. To show that (2.9) is also satisfied, we verify that $[\xi_1, \xi_2, \ldots, \xi_K]\mathbf{\Gamma} = \mathbf{B}$, i.e., that

$$\xi_1 r(M, M) + \sum_{i=2}^{K} \xi_i k_i r^{(i)}(M, M) = r(M, M+1) \tag{2.26}$$

and

$$\xi_1 k_i r^{(i)}(M, M) + \xi_i r^{(i)}(M, M) = k_i r^{(i)}(M, M+1) \text{ for } i = 2, \ldots, K \tag{2.27}$$

The equations (2.27) are trivial to verify. To verify (2.26), note that

$$\xi_1 r(M, M) + \sum_{i=2}^{K} \xi_i k_i r^{(i)}(M, M) - r(M, M+1)$$

$$= \sum_{i=2}^{K} \left\{ \left( \left(1 - \frac{M\tau^{(i)}}{\theta^{(i)}}\right) \xi_i \gamma_{i,i} - \left(1 - \frac{(1+M)\tau^{(i)}}{\theta^{(i)}}\right) \gamma_{1,i} \right.\right.$$

$$\left. + \left(1 - \frac{M\tau^{(i)}}{\theta^{(i)}}\right) \xi_1 \gamma_{1,i} \right)$$

$$\left.\times \frac{\tau^{(1)} \left(\gamma_{1,i} + M\left((1-\xi_1)\gamma_{1,i} - \xi_i \gamma_{i,i}\right)\right)^2}{\theta^{(1)} \gamma_{i,i} \left((1-\xi_1)\gamma_{1,i} - \xi_i \gamma_{i,i}\right)} \right\}$$

$$= 0$$

where the second equality holds because

$$\left(1 - \frac{M\tau^{(i)}}{\theta^{(i)}}\right) \xi_i \gamma_{i,i} - \left(1 - \frac{(1+M)\tau^{(i)}}{\theta^{(i)}}\right) \gamma_{1,i} + \left(1 - \frac{M\tau^{(i)}}{\theta^{(i)}}\right) \xi_1 \gamma_{1,i} = 0 \text{ for } i = 2, \dots, K.$$

∎

The covariance structure in (2.18)-(2.20) reflects the definition of $X(\cdot)$ as a weighted sum of independent processes. We recognize the individual summands on the left-hand side of (2.21) as squared correlation coefficients. Given (2.18)-(2.20), their sum must be less than unity if the $\gamma'_{ij}s$ are indeed covariances.

If $\gamma_{i,j} \equiv Cov\left(X^{(i)}(M), X^{(j)}(M)\right) = 0$ when $j \neq i$ and $i, j = 2, \dots, K$ as required by (2.20), and we let $X^{(1)}(\cdot) \equiv k_1^{-1}\left(X(\cdot) - \sum_{i=2}^{K} k_i X^{(i)}(\cdot)\right)$ where the $k_i's$ re given by (2.22), then we easily deduce that

$$Cov\left(X^{(1)}(M), X^{(j)}(M)\right) = 0 \text{ for } j = 2, \dots, K$$

consistently with the assumption that all the $X^{(i)}(\cdot)$'s are mutually independent.

## 3 The M-GMSI Queue

Building on the definitions of Section 1.1, we now add probabilistic assumptions. We will call the resulting model the M-GMSI queue. For some fixed $M \in \{1, 2, \dots, N-1\}$, let

$$\overline{w}_r \equiv [\overline{w}_{r,1}, \overline{w}_{r,2}, \dots, \overline{w}_{r,K}] \text{ for } r = 1, \dots, N-M \tag{3.1}$$

where

$$\overline{w}_{r,1} = \sum_{n=r}^{r+M-1} (z_{n,1} - \rho) \text{ and } \overline{w}_{r,i} = \sum_{n=r}^{r+M-1} z_{n,i} \text{ for } i = 2, \dots, K \tag{3.2}$$

Let $P^*(\cdot)$ denote the probability measure under the M-GMSI model that we will define. For $r = 1, \dots, N-M$, then let

$$H_r(z) \equiv P^*(z_{M+r,1} \leq z \mid \bar{w}_r, q_{M+r-1}) \; for \; -\infty \leq z \leq \infty, \tag{3.3}$$

$$H_r(z \mid q) \equiv P^*(z_{M+r,1} \leq z \mid \bar{w}_r, q_{M+r-1}, q_{M+r} = q) \; for \; -\infty \leq z \leq \infty, \tag{3.4}$$

$$F_r(q) \equiv P^*(q_{M+r} \leq q \mid \bar{w}_r, q_{M+r-1}) \; for \; 0 \leq q \leq \infty, \tag{3.5}$$

and

$$F_r(q \mid z) \equiv P^*(q_{M+r} \leq q \mid \bar{w}_r, q_{M+r-1}, z_{M+r,1} = z) \; for \; 0 \leq q \leq \infty. \tag{3.6}$$

By Bayes' Theorem, these distributions are related by

$$d_z H_r(z \mid q) = \frac{d_q F_r(q \mid z) d_z H_r(z)}{d_q F_r(q)} \tag{3.7}$$

for each $r$.

We will let $=_d$ denote equality of finite-dimensional distributions and assume that, in (1.1),

$$Z(t) = \rho t + \tilde{X}(t) \text{ for } t \geq 0 \text{ where } \tilde{X}(\cdot) = \tilde{X}^{(1)}(\cdot) + \sum_{i=2}^{K} k_i \tilde{X}^{(i)}(\cdot) \tag{3.8}$$

and where the $\tilde{X}^{(i)}(\cdot)'s$ are independent processes. For each $i = 1, \ldots, K$ and $r = 1, \ldots, N - M$, we will further assume that

$$\{\tilde{X}^{(i)}(t + r - 1) - \tilde{X}^{(i)}(r - 1): 0 \leq t \leq M + 1\} =_d \{X^{(i)}(t): 0 \leq t \leq M + 1\} \tag{3.9}$$

for some $(\theta^{(i)}, \tau^{(i)})$ GMSI process $X^{(i)}(\cdot)$. Implicitly $\tau^{(i)} \leq \theta^{(i)}/(M+1)$. It follows from (3.8) and (3.9) that

$$\{\tilde{X}(t + r - 1) - \tilde{X}(r - 1): 0 \leq t \leq M + 1\} =_d \{X(t): 0 \leq t \leq M + 1\} \tag{3.10}$$

where $X(\cdot)$ is a $(\theta, \tau)$ GMSI process as defined in (2.2). If, in addition,

$$z_{n,i} \equiv \tilde{X}^{(i)}(n) - \tilde{X}^{(i)}(n - 1) \text{ for } n = 1, \ldots, N \text{ and } i = 2, \ldots, K, \tag{3.11}$$

then, by (2.5), (3.3), and (3.8)-(3.11),

$$H_r(z) = N(z; \rho_r, \theta - \tau_r) \tag{3.12}$$

where

$$\rho_r = \rho - \frac{\tau^{(1)}(\bar{w}_{r,1} - \sum_{i=2}^{K} k_i \bar{w}_{r,i})}{\theta^{(1)} - M\tau^{(1)}} - \sum_{i=2}^{K} \frac{k_i \tau^{(i)} \bar{w}_{r,i}}{\theta^{(i)} - M\tau^{(i)}} \text{ and } \tau_r = \sum_{i=1}^{K} k_i^2 \frac{\theta^{(i)} \tau^{(i)}}{\theta^{(i)} - M\tau^{(i)}}.$$

Although the distribution of $z_{M+r,1}$ in (3.3) is conditional on $q_{M+r-1}$, knowing that queue length without also knowing a reference such as $q_0$ provides no information about $z_{M+r,1}$. An example where (3.9) and (3.10) are satisfied is where $\tilde{X}^{(i)} \equiv X^{(i)}$ for $i = 1, \ldots, K$ and $N \leq \min_i \theta^{(i)}/\tau^{(i)}$, as follows from the stationarity of increments of a GMSI process.

When (1.1)-(1.3) and (3.8)-(3.11) hold, the distribution in (3.6) is given by

$$F_r(q \mid z) = \begin{cases} 1 - e^{-2q(q-z)/\theta}, & q \geq q_{M+r-1} + z \text{ and } q \geq 0 \\ 0, & otherwise. \end{cases} \tag{3.13}$$

as follows from Corollary 7 of Fendick [1]. According to Proposition 3 of Fendick [1], a process $Z(\cdot)$ satisfying (3.8)-(3.10) and conditioned on a net change of $z$ over any unit interval will have the distribution over the interval of a scaled Brownian bridge with drift equal to $z$ and scaling parameter equal to $\theta$, regardless of the value of the original drift $\rho$, the parameter $\tau$, or prior history on which the process is conditioned. Since the distribution in (3.6) is conditional on $z_{M+r,1} = z$, (3.13) has the interpretation of the queue-length distribution at the end of a unit interval when the net-input process over the interval is a scaled Brownian bridge with scaling parameter $\theta$ and drift $z$ and when the queue length at the start of the interval is the known value $q_{M+r-1}$. Hajek [13] originally derived this distribution for the case in which the queue is empty at the start of the interval. From (1.5), we see that (3.13) is consistent with (1.6)

Integrating both sides of (3.7), we obtain

$$d_q F_r(q) = \int_{-\infty}^{\infty} d_q F_r(q \mid \mathfrak{z}) d_{\mathfrak{z}} H_r(\mathfrak{z}), \qquad (3.14)$$

from which it follows using (3.12) and (3.13) that

$$F_r(q) = \frac{1}{2}\left(1 - e^{\frac{-2q(\tau_r q - \theta \rho_r)}{\theta^2}} - erf\left(\frac{-q + q_{M+r-1} + \rho_r}{\sqrt{2(\theta - \tau_r)}}\right)\right. \\ \left. + e^{\frac{-2q(\tau_r q - \theta \rho_r)}{\theta^2}} erf\left(\frac{\theta(q + q_{M+r-1}) - (2\tau_r q - \theta \rho_r)}{\theta\sqrt{2(\theta - \tau_r)}}\right)\right) \qquad (3.15)$$

for $q \geq 0$ where $erf(z) \equiv 2\pi^{-1/2} \int_0^z \exp(-t^2)\, dt$. For an alternative derivation of (3.15), see Theorem 3 of Fendick [1]. Since (3.12), (3.13), and (3.15) uniquely specify the three functions on the right-hand side of (3.7), the distribution $H_r(\cdot \mid q)$ in (3.4) is uniquely specified by (3.7).

For $q \geq 0$, let

$$G_r(z \mid q) = \begin{cases} \dfrac{ae^{2q(z-q)/\theta}(4q^2 - 2qz + \theta)}{2q^2}, & z \leq q - q_{M+r-1} \\ 1, & z > q - q_{M+r-1} \end{cases} \qquad (3.16)$$

where

$$a \equiv \left(1 + \frac{e^{-2qq_{M+r-1}/\theta}(2qq_{M+r-1} + \theta)}{2q^2}\right)^{-1}. \qquad (3.17)$$

Our next proposition provides an alternative expression for $H_r(\cdot \mid q)$.

*Proposition 2: If $H_r(\cdot)$ satisfies (3.12), $F_r(\cdot \mid z)$ satisfies (3.13), $H_r(\cdot \mid q)$ is defined as in (3.4), Y is a random variable with distribution $G_r(\cdot \mid q)$ as defined in (3.16), and U is a random variable – independent of Y – with uniform distribution on $[0,1]$, then*

$$H_r(z \mid q) = P\left(Y \leq z \mid U \leq e^{-(Y-\rho_r)^2/(2(\theta-\tau_r))}\right). \quad (3.18)$$

Proof of Proposition 2:

Let $h(\cdot)$ denote any probability density function, and $H(\cdot)$ the associated cumulative distribution function. If there exists a constant $c > 0$ and probability density function $g(\cdot)$ with the same domain as $h(\cdot)$ such that

$$h(z) \leq cg(z) \quad (3.19)$$

for all $z$ in that domain and if Y has density $g(\cdot)$ and U has the uniform distribution function on $[0,1]$ where U and Y are independent, then it is well known that $H(z) = P(Y \leq z \mid U \leq h(Y)/(c \cdot g(Y)))$; see for example, page 162 of Kohlas [14]. If $H(\cdot) \equiv H_r(\cdot \mid q)$ in (3.7), and $a$ is defined by (3.17), then it follows from (3.13) that $g(z) \equiv \frac{a \cdot dF_r(q \mid z)}{dq}$ is a probability density (with respect to z) with the same domain as $h(z) = dH_r(z \mid q)/dz$. If $c \equiv$

$\left(a\sqrt{2\pi(\theta-\tau_r)}\frac{d}{dq}F_r(q)\right)^{-1}$, then, by (3.7), $h(z)/(c \cdot g(z)) =$

$e^{-(z-\rho_r)^2/(2(\theta-\tau_r))} < 1$, so that (3.19) is satisfied. ∎

The definition of the M-GMSI queuing model that follows depends on (3.8)-(3.11) only through (3.12) and (3.13). Henceforth, we will regard (3.12) and (3.13) as axioms for the M-GMSI model, for which (3.8)-(3.11) provides motivation and from which (3.15) and (3.18) follow.

## 3.1 Model Definition

Given with above definitions, the M-GMSI queuing model is defined by the following recursion:

*M-GMSI Queuing Model:*

  *Given*

- *Integers M and N such that $1 \leq M < N$*
- *Real parameters $\theta^{(i)} > 0$ and $\tau^{(i)} \leq \theta^{(i)}/(M+1)$ for $i = 1, \ldots, K$, $k_i$ for $i = 2, \ldots, K$, and $\rho$*

- Sample queue length $q_M \geq 0$
- Sample increments $z_{n,1}$ of the queue's net-input process for $n = 1, \ldots, M$
- Samples $z_{n,i}$ for $n = 1, \ldots, N-1$ and $i = 2, \ldots, K$

For $1 \leq r \leq N - M$:

- The queue length $q_{M+r}$ is a random sample from $F_r(\cdot)$ given by (3.15);
- The increment $z_{M+r,1}$ is a random sample from $H_r(\cdot \mid q_{M+r})$ given by (3.18);
- $l_{M+r} = q_{M+r} - q_{M+r-1} - z_{M+r,1}$.

If we let $\mathfrak{F}_r$ denote the sigma field generated by all known samples immediately prior to the $r^{th}$ step above, then

$$P^*(z_{M+r,1} \leq z \mid \mathfrak{F}_r) = \int_{-\infty}^{z} \int_{0}^{\infty} d_q F_r(q) \, d_{\mathfrak{z}} H_r(\mathfrak{z} \mid q) = H_r(z) \qquad (3.20)$$

$r = 1, \ldots, N - M$, where the first equality is implied by the above recursion, and the second equality by (3.7). In other words, under the M-GMSI model, the distribution of $z_{M+r,1}$ conditional on all past increments is the same as its distribution conditional on the preceding $M$. This, in addition to (3.12) and (3.13), is an axiom implicit in the model's construction. In this respect, the net-input increments of the M-GMSI model can differ from the increments of a GMSI process with drift, and this difference makes the M-GMSI model the more flexible for applications. As an illustration, the time domain of a $(\theta, \tau)$ GMSI process $X(\cdot)$ is limited when $\tau > 0$ to the interval $[0, \theta/\tau]$ by the requirement that $E[X^2(t)] = t(\theta - \tau t) \geq 0$, so that there is an upper limit to the number of samples that can be drawn at unit intervals from such a GMSI process. In contrast, regardless of the sign of $\tau$, there is no upper limit to the value $N$ for the M-GMSI model.

We note that an equivalent definition of the M-GMSI model is the recursion where, for $r = 1, ., \ldots, N - M$, $z_{M+r,1}$ is a random sample from $H_r(\cdot)$ as defined in (3.12) and $q_{M+r}$ is then a random sample from $F_r(\cdot \mid z_{M+r,1})$ as defined in (3.13). The definition we have chosen is better suited for modifications that we will introduce in Section 5 to estimate censored demands based on observed queue lengths.

We will explore further properties of the M-GMSI model through Monte Carlo simulation.

### 3.1.1 Pseudo-Random Samples of Queue Lengths

It is straightforward to generate pseudo-random samples of $q_n$ as defined by the M-GMSI model through numerical inversion of the distribution function $F_r(\cdot)$ defined by (3.15); see for example Theorem 3 on page 161 of Kohlas [14]. Note, however, that doing so requires a high-precision implementation of the error function, an example of which is algorithm 5666 from Hart et al. [15].

### 3.1.2 Pseudo-Random Samples of the Net-input Process

Proposition 2 provides the basis for the following algorithm for simulating the pseudo-random samples $z_{M+r,1}$:

*Step 1.*   *Generate Y from $G_r(\cdot \mid q_{M+r})$ in (3.16).*
*Step 2.*   *Generate U from the uniform distribution on the interval $[0,1]$*
*Step 3.*   *If $U \leq e^{-(Y-\rho_r)^2/(2(\theta-\tau_r))}$, then $z_{M+r,1} = Y$;   else return to Step 1.*

This is an example of the acceptance-rejection method for pseudo-random number generation.

## 3.2  Examples

Below, we examine simulations of queue lengths $\{q_n\}$ generated by the M-GMSI model. In Section 5 we will examine simulations of the net-input increments $\{z_{n,1}\}$ generated by the M-GMSI model.

Figures 1-3 each show queue lengths $q_n$ for $n = 0, \ldots, 2000$ simulated by the M-GMSI model for the case in which $K = 1$. The three figures were all generated using the parameters and initial values from Table 1. The figures differ in the values used for the parameter $\tau$.

Figure 1, for which $\tau = -0.1$, exhibits the largest peak queue lengths and the longest busy periods, whereas Figure 3, for which $\tau = +0.1$, exhibits the smallest and shortest.

Figure 2, for which $\tau = 0$, is an example of a discrete simulation of a queue for which the underlying net-input process is a scaled Brownian motion with drift. It exhibits peak queue lengths and busy periods that are between those in Figures 1 and 3.

In Section 5, we show that it is possible to recover the parameter $\rho$, $\theta$, and $\tau$ from sample queue lengths $q_n$ in examples similar to these. This implies that no strict subset of those parameters is sufficient to characterize the differences exhibited by Figures 1-3.

## 4  Estimation Based on Samples of Net-input Increments

Continuing to build from the notation of Sections 3, let

$$\bar{y}_r = \sum_{n=r}^{r+M}(z_{n,1} - \rho) \quad for\ r = 1,2,\ldots,N-M. \tag{4.1}$$

Also let

$$\beta_i \equiv Cov\left(\bar{y}_r, \bar{w}_{r,i}\right) \ for\ i = 1,\ldots,K \tag{4.2}$$

denote the elements of an unknown population covariance vector **B** and

$$\gamma_{i,j} \equiv Cov\left(\bar{w}_{r,i}, \bar{w}_{r,j}\right)\ for\ i,j = 1,\ldots,K \tag{4.3}$$

denote the elements of a population covariance matrix $\mathbf{\Gamma}$. Note that (4.2) and (4.3) are assumed to hold for all $r = 1,2,\ldots,N-M$.

If

$$\xi_i = \langle\mathbf{B}\mathbf{\Gamma}^{-1}\rangle_i \text{ for } i = 1,\ldots,K \tag{4.4}$$

denotes the elements of the vector $\boldsymbol{\xi}$, then it is well known that

$$\boldsymbol{\xi} = \arg\min_{\tilde{\boldsymbol{\xi}}} E\left[\left(\bar{y}_r - \tilde{\boldsymbol{\xi}}\bar{\boldsymbol{w}}_r^T\right)^2\right] \text{ for } r = 1,2,\ldots,N-M \tag{4.5}$$

so that $\boldsymbol{\xi}\boldsymbol{w}_r^T$ is the best linear predictor of $\bar{y}_r$ given $\bar{\boldsymbol{w}}_r$; see for example Section 4.1 of Whittle [16].

In the statement of the next proposition, let $E^*(\cdot)$ denote the expectation operator for the M-GMSI model.

*Proposition 3: If $\hat{\mathbf{\Gamma}}$ is a consistent estimator of $\mathbf{\Gamma}$ as defined by (4.3) and $\hat{\boldsymbol{\xi}} \equiv \hat{\mathbf{B}}\hat{\mathbf{\Gamma}}^{-1}$, where $\hat{\mathbf{B}}$ is a consistent estimator of $\mathbf{B}$ as defined by (4.2), and their elements jointly satisfy*

$$\hat{\xi}_1 < 1 + \frac{1}{M} for\ i = 1,\ldots,K, \tag{4.6}$$

$$\hat{\xi}_1 + \hat{\xi}_i \frac{\hat{\gamma}_{i,i}}{\hat{\gamma}_{1,i}} < 1 + \frac{1}{M} \text{ for } i = 2,\ldots,K \tag{4.7}$$

$$\hat{\gamma}_{i,i} > 0 \text{ for } i = 1, \ldots, K \tag{4.8}$$

$$\hat{\gamma}_{1,i} = \hat{\gamma}_{i,1} \neq 0 \text{ for } i = 2, \ldots, K \tag{4.9}$$

and

$$\hat{\gamma}_{j,i} = 0 \text{ otherwise,} \tag{4.10}$$

and if the M-GMSI model parameters are chosen to satisfy

$$\hat{k}_i = \begin{cases} \hat{\gamma}_{1,i}/\hat{\gamma}_{i,i}, & i = 2, \ldots, K \\ 1, & i = 1 \end{cases} \tag{4.11}$$

$$\hat{\theta}^{(i)} = \begin{cases} \hat{\gamma}_{i,i}\left(1 + \frac{1}{M} - \left(\hat{\xi}_1 + \hat{\xi}_i \frac{\hat{\gamma}_{i,i}}{\hat{\gamma}_{1,i}}\right)\right), & i = 2, \ldots, K \\ \left(\hat{\gamma}_{1,1} - \sum_{i=2}^{K} \hat{\gamma}_{1,i}^2/\hat{\gamma}_{i,i}\right)\left(1 + \frac{1}{M} - \hat{\xi}_1\right), & i = 1 \end{cases} \tag{4.12}$$

and

$$\hat{\tau}^{(i)} = \begin{cases} \frac{\hat{\gamma}_{i,i}}{M}\left(1 - \left(\hat{\xi}_1 + \hat{\xi}_i \frac{\hat{\gamma}_{i,i}}{\hat{\gamma}_{1,i}}\right)\right), & i = 2, \ldots, K \\ \left(\hat{\gamma}_{1,1} - \sum_{i=2}^{K} \hat{\gamma}_{1,i}^2/\hat{\gamma}_{i,i}\right)\frac{1 - \hat{\xi}_1}{M}, & i = 1 \end{cases} \tag{4.13}$$

then a consistent estimator of $\xi$ as defined in (4.5) is $\hat{\xi}$, and

$$E^*(\bar{y}_r | \bar{w}_r) = \hat{\xi}\bar{w}_r^T \text{ for } r = 1, 2, \ldots, N - M. \tag{4.14}$$

For the case in which $K = 1$, (4.11)-(4.13) reduce to $\hat{k}_1 = 1$, $\hat{\theta}^{(1)} = \hat{\gamma}_{1,1}\left(1 + M^{-1} - \hat{\xi}_1\right)$ and $\hat{\tau}^{(1)} = \hat{\gamma}_{1,1}\left(1 - \hat{\xi}_1\right)/M$.

Proof of Proposition 3: Since (4.4) holds, the property that $\hat{\xi}$ is a consistent estimator follows from the continuous mapping theorem for convergence in probability, c.f., Theorems S1 and S2 on page 102 of Goldberger [12].

If it were the case that (3.8)-(3.11) held, then it would follow from (2.5) that

$$P^*(\bar{y}_r \leq y | \bar{w}_r) = N(y; \rho_r + \bar{w}_{r,1} - \rho, \theta - \tau_r) \tag{4.15}$$

where $\rho_r$ and $\tau_{x_r}$ are defined as in (3.12). Given the parameter choices in (4.11)-(4.13), the conditional expectation in (4.14) would follow from (2.7) and Proposition 1.

It must also be true that

$$E^*(\bar{y}_r | \bar{w}_r) = \int_{-\infty}^{\infty} y \, dP^*(\bar{y}_r \leq y | \bar{w}_r); \tag{4.16}$$

c.f., page 471 of Billingsley [17]. If, as assumed for the M-GMSI model, (3.12) holds, but (3.8)-(3.11) do not necessarily hold, then, by (3.12), (4.15) continues to hold. By (4.16) the expression for the corresponding conditional expectation in (4.14) must also continue to hold. ∎

It is well known that the best linear predictor is equal to the conditional expectation when the latter is linear. This leads to the following result:

*Corollary 1: Under the conditions of Proposition 3,*

$$\hat{\xi} = \arg\min_{\tilde{\xi}} E^* \left[ (\bar{y}_r - \tilde{\xi}\bar{w}_r^T)^2 \right] \text{ for } r = 1, 2, \ldots, N - M.$$

In other words, one effect of choosing parameters for the M-GMSI model as prescribed by Proposition 3 is that the best linear predictor for $\bar{y}_r$ given $\bar{w}_r$ under its probability measure is asymptotically correct for large $N$.

At first glance, the constraint in (4.10) looks restrictive. But given $k \in [3, \ldots, K]$, suppose that (4.8)-(4.10) are satisfied for $i, j = 2, \ldots, k-1$, and let $z'_{n,k}$ for $n = 1, \ldots, N$ denote another sequence of zero-mean samples. Let $\bar{w}'_{r,k} \equiv \sum_{n=r}^{r+M-1} z'_{n,k}$ for $r = 1, 2, \ldots, N - M$. Then, if $\bar{w}_{r,i}$ is defined as in (3.2) for $i = 2, \ldots, k-1$, and $\gamma'_{i,k} \equiv Cov(\bar{w}'_{r,k}, \bar{w}_{r,i})$ does not depend on $r$, we easily deduce that

$$Cov\left( \bar{w}'_{r,k} - \sum_{i=2}^{k-1} \frac{\gamma'_{i,k}}{\gamma_{ii}} \bar{w}_{r,i}, \bar{w}_{r,j} \right) = 0 \text{ for } j = 2, \ldots, k-1$$

If we let

$$z_{n,k} = \begin{cases} z'_{n,k}, & \text{for } k = 2 \text{ and } n = 1, \ldots, N \\ z'_{n,k} - \sum_{i=2}^{k-1} \frac{\hat{\gamma}'_{k,i}}{\hat{\gamma}_{ii}} z_{n,i}, & \text{for } 3 \leq k \leq K \text{ and } n = 1, \ldots, N \end{cases} \quad (4.17)$$

where the $\hat{\gamma}'_{k,i}$'s and $\hat{\gamma}_{ii}$'s are consistent estimators of the corresponding population covariances, then it follows from (4.17) and the continuous mapping theorem that a consistent estimators of $\gamma_{j,k}$ as defined by (4.3) for $j = 2, \ldots, k-1$ is $\hat{\gamma}_{j,k} = 0$. It is also clear that knowing $(\bar{w}_{r,2}, \ldots, \bar{w}_{r,k-1}, \bar{w}_{r,k})$ is equivalent to knowing $(\bar{w}_{r,2}, \ldots, \bar{w}_{r,k-1}, \bar{w}'_{r,k})$. Given a set of zero-mean samples $z'_{n,j}$ satisfying the above assumptions for $j = 2, \ldots, K$, we can therefore derive a set of samples $z_{n,j}$ for which (4.10) is correct as a consistent estimator through the recursion defined by (4.17). This recursion is an example of Gram-Schmidt othonormalisation. The sequence $\{z_{n,k}\}$ defined by (4.17) depends implicitly on the parameter $M$ used.

Although Corollary 1 depends on samples $\{z_{n,i}\}$ for each $i$ only through their partial sums, one can start in (4.17) with sequences $\{z'_{n,i}\}$ and $\{z'_{n,j}\}$ for $i \neq j$ that are constructed by time-shifting a common sequence of samples by different amounts, so that, for example, $z'_{n,i} = z'_{n-1,j}$ for $n = 2, \ldots, N$. In this way, one can construct models fitting time series data with quite general autoregressive dependencies.

## 4.1 Estimators of Covariances

The covariances in (4.3) are closely related to the elements of Indices of Dispersion as described in Fendick & Whitt [18] and Fendick, Saksena, & Whitt [19], and the considerations discussed in Section IIIB of Fendick, Saksena, & Whitt [19] apply for their estimation. In particular, consistent estimators of the covariances in (4.3) will result from setting

$$\hat{\gamma}_{i,i} \equiv \frac{\sum_{r=1}^{N-M} \overline{w}_{r,i}\overline{w}_{r,j}}{N-M} - \left(\frac{\sum_{r=1}^{N-M} \overline{w}_{r,i}}{N-M}\right)\left(\frac{\sum_{r=1}^{N-M} \overline{w}_{r,j}}{N-M}\right) \; for \; i = 2, \ldots, K,$$

$$\hat{\gamma}_{1,i} \equiv \frac{\sum_{r=1}^{N-M} \overline{w}_{r,1}\overline{w}_{r,i}}{N-M} - \left(\frac{\sum_{r=1}^{N-M} \overline{w}_{r,1}}{N-M}\right)\left(\frac{\sum_{r=1}^{N-M} \overline{w}_{ri}}{N-M}\right)$$

$$= \frac{\sum_{r=1}^{N-M} w_{r,1}\overline{w}_{r,i}}{N-M} - \left(\frac{\sum_{r=1}^{N-M} w_{r,1}}{N-M}\right)\left(\frac{\sum_{r=1}^{N-M} \overline{w}_{ri}}{N-M}\right) \; for \; i = 2, \ldots, K,$$

and

$$\hat{\gamma}_{1,1} \equiv \frac{\sum_{r=1}^{N-M} \overline{w}_{r,1}^2}{N-M} - \left(\frac{\sum_{r=1}^{N-M} \overline{w}_{r,1}}{N-M}\right)^2 = \frac{\sum_{r=1}^{N-M} w_{r,1}^2}{N-M} - \left(\frac{\sum_{r=1}^{N-M} w_{r,1}}{N-M}\right)^2.$$

where the $\overline{w}_{r,i}$'s are defined as in (3.2) and where

$$w_{r,1} = \sum_{n=r}^{r+M-1} z_{n,1}. \tag{4.18}$$

A consistent estimator $\hat{\mathbf{B}}$ of the covariance vector $\mathbf{B}$ as defined by (4.2) is obtained analogously.

## 4.2 Estimators of the Autoregressive Coefficients and Drift

Consistently with (4.8)-(4.10), we define the *sample volatility matrix* $\hat{\Gamma}$ as the symmetric matrix with upper triangular elements

$$\hat{\mathbf{\Gamma}} = \begin{bmatrix} \hat{\gamma}_{1,1} & \hat{\gamma}_{1,2} & \hat{\gamma}_{1,3} & \cdots & \hat{\gamma}_{1,K-1} & \hat{\gamma}_{1,K} \\ & \hat{\gamma}_{2,2} & 0 & \cdots & 0 & 0 \\ & & \hat{\gamma}_{3,3} & \ddots & 0 & 0 \\ & & & \ddots & 0 & 0 \\ & & & & \hat{\gamma}_{K-1,K-1} & 0 \\ & & & & & \hat{\gamma}_{K,} \end{bmatrix}$$

where the non-zero elements are estimated as in Section 4.1. If we partition $\hat{\mathbf{\Gamma}}$ as

$$\hat{\mathbf{\Gamma}} = \begin{bmatrix} \hat{\mathbf{\Gamma}}_{1,1} & \hat{\mathbf{\Gamma}}_{1,2} \\ \hat{\mathbf{\Gamma}}_{1,2}^T & \hat{\mathbf{\Gamma}}_{2,2} \end{bmatrix}$$

with respective sizes

$$\begin{bmatrix} 1 \times 1 & 1 \times (K-1) \\ (K-1) \times 1 & (K-1) \times (K-1) \end{bmatrix}$$

then

$$\hat{\mathbf{\Gamma}}^{-1} = \begin{bmatrix} (\hat{\mathbf{\Gamma}}_{1,1} - \hat{\mathbf{\Gamma}}_{1,2}\hat{\mathbf{\Gamma}}_{2,2}^{-1}\hat{\mathbf{\Gamma}}_{1,2}^T)^{-1} & -(\hat{\mathbf{\Gamma}}_{1,1} - \hat{\mathbf{\Gamma}}_{1,2}\hat{\mathbf{\Gamma}}_{2,2}^{-1}\hat{\mathbf{\Gamma}}_{1,2}^T)^{-1}\hat{\mathbf{\Gamma}}_{1,2}\hat{\mathbf{\Gamma}}_{2,2}^{-1} \\ -\hat{\mathbf{\Gamma}}_{2,2}^{-1}\hat{\mathbf{\Gamma}}_{1,2}^T(\hat{\mathbf{\Gamma}}_{1,1} - \hat{\mathbf{\Gamma}}_{1,2}\hat{\mathbf{\Gamma}}_{2,2}^{-1}\hat{\mathbf{\Gamma}}_{1,2}^T)^{-1} & \hat{\mathbf{\Gamma}}_{2,2}^{-1} + \hat{\mathbf{\Gamma}}_{2,2}^{-1}\hat{\mathbf{\Gamma}}_{1,2}^T(\hat{\mathbf{\Gamma}}_{1,1} - \hat{\mathbf{\Gamma}}_{1,2}\hat{\mathbf{\Gamma}}_{2,2}^{-1}\hat{\mathbf{\Gamma}}_{1,2}^T)^{-1}\hat{\mathbf{\Gamma}}_{1,2}\hat{\mathbf{\Gamma}}_{2,2}^{-1} \end{bmatrix}. \quad (4.19)$$

Because $\hat{\mathbf{\Gamma}}_{2,2}$ is diagonal, all the inverse operations in (4.19) reduce to computing reciprocals. Using (4.19), we obtain $\hat{\boldsymbol{\xi}}$ directly from its definition in Proposition 3. This is an example of the approach to least-squares minimization described in Section 4.2 of Whittle [16].

For $r = 1,2,\ldots,N-M$, let $y_r = \sum_{n=r}^{r+M} z_{n,1}$, and define $w_{r,1}$ as in (4.18). Then, it is well known that

$$[\xi_1, \xi_2, \ldots, \xi_{K+1}] = \underset{[\tilde{\xi}_1, \tilde{\xi}_2, \ldots, \tilde{\xi}_{K+1}]}{\arg\min} E\left[\left(y_r - \left(\tilde{\xi}_{K+1} + [\tilde{\xi}_1, \tilde{\xi}_2, \ldots, \tilde{\xi}_K]\begin{bmatrix} w_{r,1} \\ \overline{w}_{r,2} \\ \vdots \\ \overline{w}_{r,K} \end{bmatrix}\right)\right)^2\right]$$

for $r = 1,2,\ldots,N-M$ if (4.5) is satisfied for $\xi_1, \xi_2, \ldots, \xi_K$ and if

$$\xi_{K+1} = Ey_r - \xi_1 Ew_{r,1} - \xi_2 E\overline{w}_{r,2} - \cdots - \xi_K E\overline{w}_{r,K} = \rho(M+1) - \rho M \xi_1; \quad (4.20)$$

see, for example, Section 1.1.3 on page 3 of Amemiya [10]. Consequently, a consistent estimator of $\xi_{K+1}$ is

$$\hat{\xi}_{K+1} = (N-M)^{-1} \sum_{r=1}^{N-M} \left(y_r - [\hat{\xi}_1, \hat{\xi}_2, \ldots, \hat{\xi}_K]\begin{bmatrix} w_{r,1} \\ \overline{w}_{r,2} \\ \vdots \\ \overline{w}_{r,K} \end{bmatrix}\right),$$

where $\hat{\xi}_1, \hat{\xi}_2, \ldots, \hat{\xi}_K$ are estimated as prescribed by Proposition 3. Appling (4.20) again, a consistent estimator of $\rho$ is

$$\hat{\rho} \equiv \hat{\xi}_{K+1}/(M + 1 - M\hat{\xi}_1). \tag{4.21}$$

When the $\hat{\xi}_i's$ are estimated as described above, we have not seen examples in practice for which (4.6) and (4.7) are not satisfied. One could guarantee that those constraints are satisfied using constrained least squares minimization as developed by Stark & Parker [20] and Mead & Renault [21] but the resulting estimators will then generally not be consistent ones if the constraints are binding.

### 4.3 Examples

We next consider examples in which we use the M-GMSI model with a given set of parameters to simulate samples of net-input increments and then apply Proposition 3 to obtain estimates of those parameters. All simulations have the parameters from Table 1 in common. Here and in all remaining examples in this paper, we will use the estimators of covariances and autoregressive coefficients from Section 4.1 and 4.2 for the formulas of Proposition 3.

Table 2 presents the estimated parameters from seven simulations for which $\tau = -0.1$. The bottom line of Table 2 shows that the weighted average of each of the estimated parameter values over all the runs was within 23% of the corresponding value used by the M-GMSI model to generate the samples.

Table 3 shows the corresponding results for $\tau = 0.0$, and Table 4 for $\tau = +0.1$. They illustrate that the accuracy of the estimates is better for the larger values of $\tau$. The bottom lines of Tables 3 and 4 shows that the weighted average of each of their estimates was within 10% of the parameter's true value.

Tables 2-4 each show results for five realizations in which $N = 149$ because we will later consider an application in which approximately that number of samples is available. They illustrate substantial variance of the parameter estimates in those cases.

Tables 2-4 use the same values $M = 3$ in estimating parameters as was used in simulating samples on which the estimates are based. Table 5 uses different

values $M$ in estimating parameters from samples there were generated by a single simulation of length $N = 16000$ with $M = 3$ and $\tau = +0.1$. Accurate estimates for $\rho$, $\theta$, $\tau$ are obtained when the assumed values for $M$ were less than or equal to the value $M = 3$ used by the M-GMSI simulation that generated the samples. As the assumed value of $M$ increases above the value $M$ that was used to generate the samples, the absolute size of the estimate for the parameter $\tau$ decreases towards zero.

## 5 Estimation Based on Samples of Queue Lengths

Building on the definitions from Section 3, we now consider the scenario in which sample queue lengths $q_n$ are observed, but sample increments $z_{n,1}$ of the net-input process are not.

*DM Algorithm:*

*Given*

- *Integers M and N such that $1 \leq M < N$*
- *Sample queue lengths $q_n$ for $n = 0, \ldots, N$*
- *Sample increments $z_{n,i}$ for $n = 1, \ldots, N$ and $i = 2, \ldots, K$*

*Step 0:*

- $\bar{\rho} = 0, \bar{\hat{k}}_i = 0, \bar{\hat{\theta}}^{(i)} = 0,$ and $\bar{\hat{\tau}}^{(i)} = 0$ for $i = 1, \ldots, K$.
- $\bar{\bar{l}}_n = 0$ for $n = 1, \ldots, N$
- $\bar{\bar{q}}_n = q_n$ for $n = 0, \ldots, M$ and $\bar{\bar{q}}_n = 0$ for $n = M + 1, \ldots, N$
- $z_{n,1} = q_n - q_{n-1}$ for $i = 1, \ldots, N$
- *Obtain current parameter estimates $\hat{k}_i$, $\hat{\theta}^{(i)}$ and $\hat{\tau}^{(i)}$ for $i = 1, \ldots, K$ using (4.11)-(4.13) and $\hat{\rho}$ using (4.21)*

*For $1 \leq Step \leq S$*

- *For $r = 1, \ldots, N - M$, generate $z_{M+r,1}$ as a random sample from $H_r(\cdot \mid q_{M+r})$ as given by (3.12) using current parameter estimates*
- *Update current estimates $\hat{k}_i$, $\hat{\theta}^{(i)}$ and $\hat{\tau}^{(i)}$ for $i = 1, \ldots, K$ using (4.11)-(4.13) and $\hat{\rho}$ using (4.21)*
- *For $1 \leq r \leq N - M$:*
    - $\bar{\bar{l}}_{M+r} = \bar{\bar{l}}_{M+r} + \check{l}_{M+r}/S$ where $\check{l}_{M+r} = q_{M+r} - q_{M+r-1} - z_{M+r,1}$
    - $\bar{\bar{q}}_{M+r} = \bar{\bar{q}}_{M+r} + \check{q}_{M+r}/S$ where $\check{q}_{M+r}$ is a random sample from $F_r(\cdot)$ as given by (3.15) using current parameter estimates*

- $\bar{\bar{\rho}} = \bar{\bar{\rho}} + \hat{\rho}/S, \bar{\bar{k}}_i = \bar{\bar{k}}_i + \hat{k}_i/S, \bar{\bar{\theta}}^{(i)} = \bar{\bar{\theta}}^{(i)} + \hat{\theta}^{(i)}/S$, and $\bar{\bar{\tau}}^{(i)} = \bar{\bar{\tau}}^{(i)} + \hat{\tau}^{(i)}/S$ for $i = 1, \ldots, K$.

Embedded in the DM algorithm is a modification of an M-GMSI simulation in which model-generated samples of queue lengths are replaced by observations. Increments $z_{n,1}$ for the net-input process are initially guessed, then used to estimate model parameters, new samples for the increments are then simulated based on those parameters and the observed samples of queue lengths, and the process repeated. The values assumed in Step 0 for the net-input increments $z_{n,1}$ correspond to the case in which $l_n = 0$ for $n = 1, \ldots, N$ in (1.5).

The value $\breve{l}_n$ obtained in each of Steps $1, 2, \ldots, S$ can be interpreted as a sample from the conditional distribution of the censored demand over the $n^{th}$ interval given the observed queue lengths at the start and end of that interval. The average value $\bar{\breve{l}}_n$ obtained at the end of all the steps is therefore an estimate of the conditional expectation of $\breve{l}_n$. Similarly, $\breve{q}_n$ is a sample from the conditional distribution for the queue length at the end of the $n^{th}$ interval given the observed queue length at the start of that interval, and $\bar{\breve{q}}_n$ is an estimate of the conditional expectation of $\breve{q}_n$.

The samples $\{z_{n,i}\}$ for $i = 2, \ldots, K$ have the interpretation as mutually orthonormalized increments of other processes that may be predictive of the value of the unobserved net-input process. The conditional distributions and conditional expectations discussed above also depend on them when they are available.

The examples that follow use samples of queue lengths from M-GMSI simulations with the values from Table 1 in common. In particular, $K = 1$ for these examples, so that $\hat{\theta} \equiv \hat{\theta}^{(1)}$, $\hat{\tau} \equiv \hat{\tau}^{(1)}$, $\bar{\bar{\theta}} \equiv \bar{\bar{\theta}}^{(1)}$ and $\bar{\bar{\tau}} = \bar{\bar{\tau}}^{(1)}$.

## 5.1 Examples: Estimating Parameters

For sample queue lengths $q_n$ for $n = 1, 2, \ldots, 8000$ simulated using $\tau = +0.1$, Table 6 shows estimates $\hat{\rho}, \hat{\theta}$ *and* $\hat{\tau}$ obtained from the DM algorithm after Step 0 and after each of the next six steps. The challenge in estimating parameters from queue lengths is highlighted by the estimate of $\hat{\rho} = 0.001$ after Step 0. An estimate near zero for $\hat{\rho}$ is to be expected at the end of Step 0, since the average queue-length increment must be near zero if queue-lengths remain stable over the course of the simulation. Nevertheless, after each of Steps 1-6, the estimate in Table 6 for each parameter is reasonably close to its population value. The average values $\bar{\hat{\rho}}, \bar{\hat{\theta}},$ and $\bar{\hat{\tau}}$ over six steps are each within 10% of the corresponding population value.

Table 7 shows the average parameter estimates $\bar{\hat{\rho}}, \bar{\hat{\theta}},$ and $\bar{\hat{\tau}}$ obtained using the DM model after the given number of steps from samples $q_n$ simulated using varying values of $N$ and $\tau$. The accuracy is roughly the same as we obtained in Tables 2-4 based on samples of the net-input increments $z_{n,1}$.

## 5.2 Examples: Estimating Conditional Expectations

We next use the DM algorithm to produce the estimated conditional expectation function $\{\bar{\bar{l}}_n\}$ of the censored demand using samples of queue lengths obtained from an M-GMSI simulation, and we compare those conditional expectations to the sample path $\{l_n\}$ from the same M-GMSI simulation. The M-GMSI simulation here used the parameters $N = 2000$, $\tau = +0.1$. The DM model used 100 steps to obtain its estimates. Figure 4 plots the samples $l_n$ from the M-GMSI simulation above the horizontal axis and the estimates $\bar{\bar{l}}_n$ of conditional expectations below the horizontal axis. As one would expect, the sample path is more volatile than the estimated conditional expectation function. The average of the conditional expectations was within 15% of the average of the sample-path values.

Figure 5 shows the corresponding results for five realization for which $N = 149$. Again the sample path is more volatile than the estimated conditional expectation. But although the peaks of the latter are lower, the average value of the latter is greater than that of the former in one of the four cases and within 6% in two of the four cases. (It was within 45% in the other two cases.) Figure 5 illustrates that the estimated conditional expectations of censored demands are useful tools for understanding the timing and rough scale of the unobserved censored demand even when the number of samples of queue length samples is limited. Estimated quantiles of censored demand would yield further insight into sample-path behavior but are beyond the scope of this paper.

The results in Figures 4 and 5 and Tables 6 and 7 assumed that the same value for the parameter $M$ is used for the DM algorithm as was originally used for the M-GMSI simulation that generated the sample queue lengths. If that original value is not known, it is natural to try different values with the goal of using the one that maximizes the goodness of fit for the model. We will quantify the goodness of fit for the DM algorithm by comparing the estimates $\bar{\bar{q}}_n$ of the conditional expectations for queue lengths with the actual queue lengths $q_n$ for $n = M + 1, \ldots, N$. In particular, we will define the coefficient of determination as

$$R^2 \equiv 1 - \frac{\sum_{n=M+1}^{N}(q_n - \bar{\bar{q}}_n)^2}{\sum_{n=M+1}^{N}(q_n - \bar{q})^2} \quad where \quad \bar{q} \equiv \frac{\sum_{n=M+1}^{N} q_n}{N + M}. \tag{5.1}$$

The coefficient of determination describes the proportion of the variation in the $q_n$'s that is explained by the estimated conditional expectations.

Table 8 shows $\bar{\bar{\rho}}, \bar{\bar{\theta}}, \bar{\bar{\tau}}$, and $R^2$ obtained after 50 steps of the DM algorithm using different guesses for $M$ given queue lengths generated by the M-GMSI model under the same assumptions as in Table 6. (The M-GMSI simulation used $M = 3$ in this case, as in all the other examples.) Curiously, $R^2 = 0.88$ (to two significant digits) for all the values of $M$ considered in Table 8. (An example in which $R^2$ does vary with $M$ is described in Section 7.)

Table 8 tells the same story about estimates of the parameters $\rho$, $\theta$ and $\tau$ obtained using samples of queue length as Table 5 did about estimates obtained using samples of net-input increments. The estimates obtained for $M$ less or equal to the correct value ($M = 3$) roughly agree with the estimates obtained using the correct value, but estimated obtained for $M$ greater than the correct value do not. A rule of thumb for correctly choosing $M$ in this case where $R^2$ does not provide differentiation is to choose the largest value for which parameter estimates for $\rho$, $\theta$ and $\tau$ are still representative of those obtained for smaller values.

The insensitivity of $R^2$ to the parameter $M$ for the example in Table 8 is an opportunity to explore the question of whether different models that achieve the same goodness of fit of sample queue lengths (as measured by $R^2$) will produce similar estimates $\bar{\bar{l}}_n$ for the conditional expectation of censored demand in each of the $N$ periods. In Figures 6-8, we present estimates $\bar{\bar{l}}_n$ obtained from the DM algorithm given a sample path $\{q_n\}$ of queue lengths from an M-GMSI simulation that used parameters $N = 2000, \tau = +0.1$ and the common parameters from Table 1. Figure 6 compares the estimates $\bar{\bar{l}}_n$ obtained when the M-GMSI model assumed the correct value $M = 3$ to those obtained when it assumed $M = 1$. The estimates differ somewhat, as we would expect given the finite number of steps of the M-GMSI model, but look similar. Figure 7 presents the corresponding comparison for $M = 3$ vs. $M = 10$, and Figure 8 for $M = 3$ vs. $M = 50$. There is a noticeable difference in the scale of estimates of censored demand in Figure 7, which becomes more pronounced in Figure 8. Nevertheless, the timing and rough scale of censored demand is the same in all three examples.

## 6  Queues with Graduated Boundaries

The results and methods of Sections 3-5 generalize for queues with graduated boundaries. Below, we construct a family of models, indexed by a real parameter $\omega \geq 1$, assuming that (1.1), (1.2), and (3.8)-(3.11) hold, but not necessarily (1.3). Noting that (1.5) and (1.6) continue to hold under these weaker assumptions, we denote the distribution defined in (3.6) for model $\omega$ by $F_r(q \mid z)_\omega$ and assume that $F_r(q \mid z)_{\omega=1}$ corresponds to the case in (3.13) when (1.3) also holds. For a given value $q_{n-1}$ in (1.5), the censored demand $l_n$ and queue length $q_n$ in (1.5) for model $\omega$ are each at least as great as it would be in case (1.3) also held, as follows from Remark (7) on page 20 of Harrison [3]. We conclude that, for any $z$,

$$F_r(q \mid z)_\omega \leq F_r(q \mid z)_1 \equiv \begin{cases} 1 - e^{-2q(q-z)/\theta}, & q \geq q_{M+r-1} + z \text{ and } q \geq 0 \\ 0, & otherwise. \end{cases} \quad (6.1)$$

This is an example of first-order stochastic dominance as defined on page 136 of Wolfstetter [22].

A distribution function satisfying (6.1) is

$$F_r(q \mid z) = F_r(q \mid z)_\omega \equiv \begin{cases} 1 - e^{-2q(q-z)/(\omega\theta)}, & q \geq q_{M+r-1} + z \text{ and } q \geq 0 \\ 0, & otherwise. \end{cases} \quad (6.2)$$

for $\omega \geq 1$. When (6.2) holds, $F_r(q \mid z)_{\omega_2} \leq F_r(q \mid z)_{\omega_1}$ if $\omega_2 > \omega_1$, and the inequality is strict at all values $q$ for which the distributions are non-zero. Since the derivation of (3.12) did not depend on (1.3), we will continue to assume that (3.12) holds, but will replace (3.13) with (6.2). Notice that ω appears in (6.2), but not in (3.12). When (3.12) and (6.2) hold, (3.15) must be replaced by

$$F_r(q) = \frac{1}{2}\left(1 - e^{\frac{-2q(q(\tau_r+(\omega-1)\theta)-\theta\omega\rho_r)}{\omega^2\theta^2}} - erf\left(\frac{-q+q_{M+r-1}+\rho_r}{\sqrt{2(\theta-\tau_r)}}\right)\right.$$

$$\left. + e^{\frac{-2q(q(\tau_r+(\omega-1)\theta)-\theta\omega\rho_r)}{\omega^2\theta^2}} erf\left(\frac{-2\tau_r q + q\theta(2-\omega) + \omega\theta(\rho_r + q_{M+r-1})}{\omega\theta\sqrt{2(\theta-\tau_r)}}\right)\right).$$

as we derived from (3.14) with the help of Mathematica. The expression in (3.18) for $H_r(\cdot \mid q)$ will then continue to hold if we replace $\theta$ by $\omega\theta$ in (3.16) and (3.17) so that

$$G_r(z \mid q) = \begin{cases} \dfrac{ae^{2q(z-q)/(\omega\theta)}(4q^2 - 2qz + \omega\theta)}{2q^2}, & z \leq q - q_{M+r-1} \\ 1, & z > q - q_{M+r-1} \end{cases}$$

and

$$a \equiv \left(1 + \frac{e^{-2qq_{M+r-1}/(\omega\theta)}(2qq_{M+r-1} + \omega\theta)}{2q^2}\right)^{-1}.$$

With these changes, the recursion in Section 3.1 defines a generalization of the M-GMSI model, which we will call the ωM-GMSI queue. An ωM-GMSI queue has a sharp boundary if $\omega = 1$ and a graduated boundary otherwise. As $\omega$ increases, more demand is censored away from the boundary at zero.

Under the ωM-GMSI model, the random samples of censored demands $l_n$ in (1.5) will depend on ω, but the random samples of net input increments $z_{n,1}$ will not. If $\bar{l} \equiv \sum_{i=1}^{N} l_n/N$ and $\bar{z} \equiv \sum_{i=1}^{N} z_{n,1}/N$ and if queue lengths remain stable over the course of a realization, then (1.5) implies that $\bar{l} \approx -\bar{z}$ for large $n$. Therefore, the choice of ω will affect the timing of censored demands, but not so much the total amount over long intervals. If more demand is censored away from the lower bound at zero, then less demand is censored at zero to maintain non-negative queue lengths.

We gain further insight into the behavior of the ωM-GMSI model by noting that it can be equivalently defined through the recursion where, for $r = 1,.,...,N-M$, $z_{M+r,1}$ is a random sample from $H_r(\cdot)$ as defined in (3.12) and where $q_{M+r}$ is then a random sample from $F_r(\cdot \mid z_{M+r,1})$ as defined in (6.2). By (6.2) and Theorem 3 on page 161 of Kohlas [14], we can simulate samples from the latter distribution by letting

$$q_{M+r} = \begin{cases} q_{M+r-1} + z_{M+r,1}, & U \leq 1 - e^{-2q_{M+r-1}(q_{M+r-1}-z_{M+r,1})/(\omega\theta)} \\ \frac{1}{2}\left(z_{M+r,1} + \sqrt{z_{M+r}^2 - 2\omega\theta \log(1-U)}\right), & otherwise. \end{cases} \quad (6.3)$$

if $q_{M+r-1} + z_{M+r,1} \geq 0$, where $U$ is a uniformly distributed random variable on [0,1]. This shows that the choice of $\omega$ affects the sample path of queue lengths primarily within some neighborhood of the lower bound at zero. By (1.5) and (6.3), the probability that $l_n = 0$ approaches one as $q_{M+r-1}$ become large. Therefore, the queue-length process under the $\omega$M-GMSI model will still tend to alternate between intervals in which demand is being censored to some degree in successive periods and intervals in which it is not censored at all.

A theoretical question is whether, for any given distribution $F_r(q \mid z)_\omega$ satisfying (6.1), there always exists a model satisfying (1.1), (1.2) and (3.8)-(3.11) for which $F_r(q \mid z)$ as defined by (3.6) equals $F_r(q \mid z)_\omega$ for $r = 1, 2, \ldots, M - N$. Borrowing from results on page 138 of Wolfstetter [22] for general stochastically ordered random variable, a sketch of an existence proof is as follows: if $q_{M+r}^-$ is random sample from $F_r(\cdot \mid z_{M+r,1})_1$ as defined in (6.1), then it is easily verified that $q_{M+r} \equiv F_r^{-1}\left(F_r(q_{M+r}^- \mid z_{M+r,1})_1\right)_\omega$ has distribution function $F_r(\cdot \mid z_{M+r,1})_\omega$ and, by (6.1), that

$$P^*(q_{M+r} \geq q_{M+r}^- \mid z_{M+r,1}) = P^*\left(F_r^{-1}\left(F_r(q_{M+r}^- \mid z_{M+r,1})_1\right)_\omega \geq q_{M+r}^-\right)$$
$$= P^*\left(F_r(q_{M+r}^- \mid z_{M+r,1})_1 \geq F_r(q_{M+r}^- \mid z_{M+r,1})_\omega\right) = 1,$$

so that $\delta l_{M+r} \equiv q_{M+r} - q_{M+r}^- \geq 0$ with probability one. Since $F_r(\cdot \mid z_{M+r,1})_1$ is the distribution from (3.13) when $z = z_{M+r,1}$, we can interpret $q_{M+r}^-$ as the left-hand limit of the queue length as time approaches $M + r$ for a realization satisfying (1.1)-(1.3) and (3.8)-(3.11) on $(M + r - 1, M + r)$. And we can interpret $\delta l_{M+r}$ as an impulse that, if added to the cumulative censored demand at time $M + r$, results in a queue length $q_{M+r}$ with the given distribution function. The impulse is a random variable that is dependent on $q_{M+r}^-$ and $z_{M+r,1}$ and is non-negative with probability one. Although an impulse of censored demand is an abstract concept, it can be interpreted as demand that would have offset a simultaneous impulse in supply of the same size under the M-GMSI model. Strictly speaking, the assumption underlying the M-GMSI model that the net-input process is continuous does not require that the cumulative supply and demand processes are themselves continuous but only that any impulse in supply is simultaneously offset by an impulse in demand of the same size.

We can modify the DM algorithm from Section 5 also to account for graduated boundaries by substituting the functions $H_r(\cdot \mid q_{M+r})$ and $F_r(\cdot)$ in the original DM algorithm of Section 5 for the ones developed above for the ωM-GMSI model. Henceforth, we call the result of those substitutions the ω*DM algorithm*. When we refer to the *DM algorithm* without qualification, we will mean the original algorithm from Section 5. Embedded in the ωDM algorithm is a modification to an ωM-GMSI simulation in which model-generated samples of queue lengths are replaced by observations.

The choice in (6.2) is not the only one satisfying (6.1), and other choices could conceivably results in a better fit for a given data set. The above development of the ωM-GMSI model provides a template for constructing alternative models based on other choices.

## 6.1 Examples: Queue Simulations

In this section, we examine simulations of queue lengths $\{q_n\}$ for the ωM-GMSI model.

Figures 9 and 10 each show queue lengths $q_n$ for $n = 0, \ldots, 2000$ simulated by the ωM-GMSI model for the case in which $K = 1$. The two figures were both generated with parameters and initial values from Table 1 and with $\tau = 0.1$. The figures differ in the values used for the parameter $\omega$. Since the ωM-GMSI model is equivalent to the ωM-GMSI model when $\omega = 1$, Figure 3 provides the corresponding results for $\omega = 1$.

The queue lengths in Figures 9, for which $\omega = 5.0$, approach zero much less frequently than those of Figure 3, for which $\omega = 1.0$. And the queue lengths in Figure 10, for which $\omega = 10.0$, even less so. Consequently, Figures 9 and 10 illustrate how queues typically stabilize at a higher level when the boundary is graduated and how the amount of graduation increases with $\omega$. The queue lengths come close to zero at a few different points in Figures 9 and 10, but typically remain at values significantly above zero. The largest queue lengths in Figures 3, 9, and 10 have roughly the same scale. That is consistent with earlier observations that queues with different values $\omega$ evolve similarly at a sufficient distance from zero.

Figure 11 shows censored demands $l_n$ and queue lengths $q_n$ generated by the same $\omega$M-GMSI simulation for the case in which $\omega = 1.0$ and in which the other parameters are the same as in Figures 9 and 10. This is an example of a queue with a sharp boundary, and Figure 11 shows that its demand is censored only when nearby queue lengths are near zero. Figure 12 shows corresponding results for $\omega = 5.0$, and Figure 13 for $\omega = 10$. Demand is typically censored farther from the boundary in Figure 12 than in Figure 11, and farther still from the boundary in Figure 13.

In experiments with long running simulation (e.g., $N = 8000$), we have verified that the average of censored demands (over all periods) is roughly equal to $-\rho$ regardless of $\omega$. Because of the relatively short durations of the simulations in Figures 11, 12, and 13, the average value of censored demand differs substantially for the three cases.

## 6.2 Examples: Estimating Parameters

Table 9 shows estimates obtained after 50 steps of the $\omega$DM algorithm using different guesses for $\omega$ given queue lengths generated by an $\omega$M-GMSI simulation with $N = 8000, M = 3, \omega = 5.0, \rho = -0.1, \theta = 1.1,$ and $\tau = +0.1$. For the example in Table 9, the $\omega$DM algorithm produced accurate estimates for $\rho, \theta,$ and $\tau$ when it used the correct value of $\omega = 5.0$, but not otherwise.

The $R^2$ value in Table 9 was computed based on the formula in (5.1). Since $R^2 = 0.85$ (to two significant digits) for all the values of $\omega$ considered in Table 9, the $R^2$ value did not help in identifying which value of $\omega$ was the correct one.

Table 9 also shows "STD % Error", which was computed as follows: for each row in Table 9, queue lengths were simulated using the $\omega$M-GMSI model with $N = 8000, M = 3$, and the row's values $\omega, \rho = \overline{\rho}, , \theta = \overline{\overline{\theta}}$ and $\tau = \overline{\tau}$. The standard deviation of those queue lengths over the length of the simulation run was next calculated and compared to the standard deviation of the observed queue lengths $q_n$ over all $n$; and the % error reported as "STD % Error". Consequently, STD % Error compares a global property of the model for each value $\omega$ to a global property of the observations. The minimum value for STD % Error in Table 9 occurred at the value $\omega = 5$ originally used in generating the observed queue lengths.

# 7  An Application

We will now apply the methods of this paper to study the problem of unfulfilled demand for jobs using publicly available data from the U.S. Department of Labor: Bureau of Labor Statistics. In particular, we study total U.S. non-farm job openings, job separations, and hires, not seasonally adjusted, as reported for each of the 151 months from December 2000 through June 2013. For month $n = 1, \ldots, 151$, $p_n$ will denote the reported job openings at the end of the month (a stock variable), $s_n$ the reported job separations over the month (a stock variable), and $h_n$ the reported hires over the month (another stock variable). According to the National Bureau of Economic Research, recessions occurred from March 2001 through November 2001 and from December 2007 through June 2009. To aid in interpreting the data, we will show the month name rather than the month index in the figures that follow.

Our goal is to estimate the censored demand for jobs each month, i.e., to estimate how many job openings would have been filled each month had there been more of the same openings. To apply the results developed in this paper to solve this problem, we must create a model in which the number of openings corresponds to a queue length, and changes in that queue length are driven by a net-input process representing the difference between supply and demand for job openings. The demand in this case is the demand to close job openings and is comprised of the demand of job seekers to fill jobs and the demand of job advertisers to close them without having filled them (e.g., because of changing circumstances within a firm). It is not meaningful to speak of unsatisfied demand of the later type, so we will attribute all censored demand to the former type.

For each month, Figure 14 shows the openings $p_n$, Figure 15 the separations $s_n$, and Figure 16 the hires $h_n$. The separations and hires (particularly the later) exhibit seasonality, a point to which we return later. Although the Bureau of Labor Statistics also provides seasonally adjusted data, we consider here only non-seasonally adjusted values to remain closer to the raw data and to highlight the flexibility of our models.

Figure 14 shows that the number of openings have remained significantly above zero over the 13-year period, even during recessions. We gain further insight by plotting the ratio $p_n/h_n$ of openings to hires in Figure 17. This ratio can be interpreted as an estimate of the average time to fill a job opening at the end of month $n$. There must exist an upper limit on how quickly, on average, persons looking for a job can identify a suitable one and then complete the interview and hiring processes. Or, equivalently, there must exist a lower bound imposed by logistics on the average number of months to fill a job opening. Figure 17 suggests that during certain periods, the average time to fill a job opening approached half a month, which would seem to be close to a lower bound on what is possible. Since there is no reason to expect the lower bound to remain constant over long periods, Figure 17 also plots the line $\alpha_0 + \beta_0 n$, where $\alpha_0$ and $\beta_0$ minimize $\sum_n (p_n/h_n - (\alpha_0 + \beta_0 n))^2$ subject to the constraint that $\alpha_0 + \beta_0 n \leq p_n/h_n$ for all $n$. This line is our estimate of the lower limit on the time to fill an opening as a function of $n$ (month). Given the number of hires $h_n$ reported for month $n$, the quantity $(\alpha_0 + \beta_0 n)h_n$ is an estimate of the expected number of openings at the end of month $n$ created more recently than the lower bound on the expected time to fill an opening. Consequently, one should expect to see at least $(\alpha_0 + \beta_0 n)h_n$ openings at the end of month $n$ regardless of the demand for jobs at that time. To derive a queue length process that conforms as closely as possible to the assumptions of Section 1.2, we therefore define the *fillable* openings as

$$q'_n \equiv p_{n+1} - (\alpha_0 + \beta_0(n+1))h_{n+1} \geq 0 \ for \ n = 0,2,...N,$$

where $N = 150$ in this case. Henceforth, we will work with the the fillable openings rather than with the total openings. Figure 18 plots the fillable job openings, as well as predictions for them that we discuss later.

When the unobserved net-input process has a constant drift and stationary covariance structure, we would expect the queue length process also to appear stationary. In practice, it is reasonable to expect that both the supply and demand for jobs will increase exponentially over long time scales as determined by growth of the economy and the population, so that there may be a subtle exponential trend in the fillable openings $\{q'_n\}$. To correct for it, we let

$$q_n \equiv \frac{q'_n}{exp(\alpha_1 + \beta_1 n)} \ for \ n = 0, 2, \ldots, N \tag{7.1}$$

where $\alpha_1$ and $\beta_1$ minimize $\sum_{n=0}^{N}(\ln q'_n - (\alpha_1 + \beta_1 n))^2$. With this transformation, $q_n$ remains non-negative, and $q_n = 0$ if and only if $q'_n = 0$. Although $\{s_n\}$ and $\{h_n\}$ are flow variables, they are also defined to be non-negative; so we will similarly de-trend and center them by letting

$$z'_{n,2} \equiv \frac{s_n}{exp(\alpha_2 + \beta_2 n)} - 1 \ and \ z'_{n,3} \equiv \frac{h_n}{exp(\alpha_3 + \beta_3 n)} - 1 \ for \ n = 1, \ldots, N \tag{7.2}$$

where $\alpha_2$ and $\beta_2$ minimize $\sum_{n=1}^{N}(\ln s_n - (\alpha_2 + \beta_2 n))^2$ and $\alpha_3$ and $\beta_3$ minimize $\sum_{n=1}^{N}(\ln h_n - (\alpha_3 + \beta_3 n))^2$. For given $M$, we finally obtain the samples $\{z_{n,2}\}$ and $\{z_{n,2}\}$ by applying (4.17) to orthonormalize the values obtained from (7.2).

Through the above steps, we transformed the raw data into suitable input for the DM algorithm. To interpret the output of the DM algorithm, we must then invert the transformation in (7.1). We do so by letting

$$\bar{q}'_n = \bar{q}_n exp(\alpha_1 + \beta_1 n) \ for \ n = 0, \ldots, N, \tag{7.3}$$

$$\breve{l}'_n = \breve{l}_n \left( \frac{exp(\alpha_1 + \beta_1(n-1)) + exp(\alpha_1 + \beta_1 n)}{2} \right) \ for \ n = 1, \ldots, N, \tag{7.4}$$

and

$$\bar{l}'_n = \bar{l}_n \left( \frac{exp(\alpha_1 + \beta_1(n-1)) + exp(\alpha_1 + \beta_1 n)}{2} \right) \ for \ n = 1, \ldots, N. \tag{7.5}$$

The DM algorithm in Section 5 defines the terms that are rescaled on the right-hand side of (7.3)-(7.5). The DM algorithm produces samples $\breve{l}_n$ and estimates $\bar{l}_n$ for each $n$ conditionally on the observed values $q_{n-1}$ and $q_n$. If it were the case that $\beta_1 = 0$, so that the queue lengths were originally scaled by the same amount in (7.1), then $\breve{l}_n$ and $\bar{l}_n$ also would be scaled by that amount; and the rescaling used in (7.4) and (7.5) would be exact. More generally, we would expect $exp(\alpha_1 + \beta_1(n-1))$ and $exp(\alpha_1 + \beta_1 n)$ to be close to one another for each $n$ in applications, so that (7.4) and (7.5) should be reasonable approximations.

Using the definition in (7.3), we will define the coefficient of determination as

$$R^2 \equiv 1 - \frac{\sum_{n=M+1}^{N}(q'_n - \bar{\bar{q}}'_n)^2}{\sum_{n=M+1}^{N}(q'_n - \bar{q}')^2} \quad where \quad \bar{q}' \equiv \frac{\sum_{n=M+1}^{N} q'_n}{N-M}. \tag{7.6}$$

Table 10 shows estimates obtained from the DM algorithm after 500 steps for the case in which $K = 1$, so that the input comprised only the samples $\{q_n\}$. Table 11 shows the corresponding estimates for the case in which $K = 3$ when the input comprises the samples $\{q_n\}, \{z_{n,2}\}$ and $\{z_{n,2}\}$. We will refer to the case in which $K = 1$ as the one-dimensional fit, and the case in which $K = 3$ as the three-dimensional fit. In the later case, the DM algorithm produces estimates $\bar{\bar{k}}_i, \bar{\theta}^{(i)}$ and $\bar{\tau}^{(i)}$ for $i = 1,2,3$, and we obtain the values presented in Table 11 by letting

$$\bar{\bar{\theta}} \equiv \sum_{i=1}^{K} \bar{\bar{k}}_i^2 \bar{\theta}^{(i)} \quad and \quad \bar{\bar{\tau}} \equiv \sum_{i=1}^{K} \bar{\bar{k}}_i^2 \bar{\tau}^{(i)},$$

as guided by to (2.3).

Figure 19 plots the $R^2$ values from Tables 10 and 11. There is close agreement between the one- and three dimensional models except for the cases $M = 9, 10,$ and $11$. For reference, a naïve model for which $\bar{\bar{q}}'_n \equiv q'_{n-1}$ in (7.6) results in $R^2 = 0.6$. The estimates in Figure 19 are close to that value in some cases. The case in Figure 19 in which $M = 1$ is the only one for which the $R^2$ estimate from the three-dimensional model is below that from the one-dimensional model, and those two estimates are close to one another. We expect that the three-dimensional fit would be uniformly better if a sufficient number of steps were used by the DM algorithm.

For the one-dimensional model, the choice $M = 1$ results in a marginally better fit than the choice $M = 11$. We conclude that $M = 1$ provides the best fit among the one-dimensional cases.

For the cases in Figure 19 for which $M = 9, 10,$ and $11$, the $R^2$ values obtained with the three-dimensional fit are substantially greater than those obtained from the one-dimensional model, with the highest value of $R^2 = 0.78$ occurring for the three dimensional fit when $M = 11$ months. As with any regression analysis, there is a risk that adding dependent variables to the model will introduce spurious dependencies. As an example, the $R^2$ values of 0.88 reported in Table 8 must be close to the maximum that can be meaningfully obtained in those cases, since it is the value obtained using the same dimension $K = 1$ as was originally used in simulating the samples. Adding dependent variables to those models would necessarily introduce spurious dependences. But given the seasonality of separations and hires discussed earlier, it is understandable why conditioning on the last 11 months of data for the three-dimensional model of fillable job openings might provide substantial information about what will happen the next month. Therefore, the fact that the maximum $R^2$ value for the three-dimensional model occurs at $M = 11$ suggests that the dependence on the variables $\{z_{n,2}\}$ and $\{z_{n,2}\}$ is meaningful and that that the best-fit case for the three-dimensional model occurring at $M = 11$ represents a meaningful improvement relative to the best-fit case for the one-dimensional model.

For the best-fit three-dimensional model, Figure 18 shows for each month $n$ the estimated conditional expectation $\bar{\bar{q}}'_n$ of the fillable openings given the previously observed fillable openings $q'_{n-1}$ and previously observed samples of separations and hires. The estimates lag behind changes in the fillable openings in some periods but anticipate them in others, including the recessionary period in the first half of 2009 when the fillable openings dropped sharply. Because of the way that way dates were mapped months to indices for job openings, separations, and hires, our predictions for openings each month in Figure 18 depended on separations and hires two months prior, but not one month prior. For example, the fillable openings predicted for the end of March 2004 depended on separations and hires from January 2004 (as well as the ten months before that) but not from February 2004. This resulted in a slightly better fit than the natural alternative of predictions using one-month prior separations and hires.

Figure 20 shows the estimated conditional expectations $\bar{\bar{l}}'_n$ from (7.5) for the censored demand for jobs each month obtained using the best-fit one-dimensional and three-dimensional models described above. Since, as we just observed, the best-fit three-dimensional model is substantially better, we take its results for the estimates $\bar{\bar{l}}'_n$ to be the more accurate.

Figure 20 can be applied to answer the narrow question whether, given prior history, more of the same openings would have resulted in more hires in that month. The results there would support a hypothesis that a greater supply of such jobs would have resulted in significantly more immediate hires during 2009 and mid-2010 (when the unemployment rate remained in the neighborhood of 10%) and, to a lesser extent, during periods of 2002 and 2003. The results there do not support a hypothesis that more of the same jobs would have resulted in more immediate hires in other periods. Since the results in Figure 20 for each month are conditional on the actual history over prior months, they do not directly address the question of whether a greater supply of openings in any month would have had a longer-term impact on supply or demand for jobs.

In Figure 20, the estimated expected censored demands are substantially greater for the best-fit three-dimensional model than for the best-fit one-dimensional model. We have also observed in additional examples that the estimates of expected censored demands tend to increase with the $R^2$ value. As a result, a conservative approach to including dependent variables in the model will tend to result in estimates of censored demands that are lower bounds on the true values. Conversely, a model that artificially inflates $R^2$ by arbitrarily including dependent variables is likely to overestimate the true censored demand. Since we have not investigated adding more dependent variables to the model, we believe that there is some room for meaningful improvement to $R^2$ and that such an improvement would likely result in a tighter lower bound on expected censored demands. As discussed in reference to Figures 4 and 5, we would expect that the estimated expected censored demand in Figure 20 is less volatile than the actual (unobserved) sample path of censored demand.

The final question that we address in this paper is whether a model with a graduated boundary would result in different estimates of expected censored demand and might better explain the data. Since we are interested in understanding the sensitivity of the preceding results to the assumption of a sharp boundary, we use a three-dimensional model with the same dependent variables as used above and the best-fit parameter of $M = 11$ obtained from the original DM algorithm.

Figure 21 compares estimates of the expected censored demand $\bar{\bar{l}}'_n$ obtained after 1000 steps of the ωDM algorithm with $\omega = 1$ and $\omega = 2$, and Figure 22 shows the corresponding results with $\omega = 1$ and $\omega = 4$. The two cases in Figure 21 look similar, though the differences are noticeable. The differences are more pronounced in Figure 22, where the case of $\omega = 4$ shows censoring of demand over a broader range of dates (including most of 2011) than does the case with $\omega = 1$.

Table 12 shows estimates for $\rho$, $\bar{\bar{\theta}}$, and $\bar{\bar{\tau}}$ (defined in the same way as in Table 11) and for $R^2$ (defined as in (7.6)) after 500 steps of the ωDM algorithm with varying values ω. The other assumptions are the same as in Figures 21 and 22. The case for which $\omega = 1$ is equivalent to the original DM algorithm, and the small differences reported for that case in Table 12 from the case in Table 11 for which $M = 11$ are easy attributable to estimation error.

Table 12 shows that the $R^2$ value achieves its maximum at $\omega = 1$ and decreases as ω increases. On the other hand, $R^2$ is relatively flat over the region $1 \leq \omega \leq 2$ and therefore does not strongly distinguish between models with ω in that range. The $R^2$ value for $\omega = 4$ is substantially lower than for $\omega = 1$.

Table 12 also shows results for the standard deviation of queue lengths ("STD"). As in Table 9, queue lengths were simulated using the ωM-GMSI model but, in the cases here, assuming the availability of samples of separations $z_{n,2}$ and hires $z_{n,3}$ for $n = 1, \ldots 150$. Since the availability of those sample limit each realization to a duration of $N = 150$, we generate ten independent realizations for each value ω. To create each realization, the estimates $\bar{\bar{k}}_i$, $\bar{\bar{\theta}}^{(i)}$ and $\bar{\bar{\tau}}^{(i)}$ for $i = 1,2,3$ are obtained from the observed queue lengths $q_n$ for $n = 0, \ldots, 150$ (the fillable job openings) and the other available samples using 500 steps of the ωDM algorithm with $M = 11$. The ωM-GMSI model then uses those parameters (including the parameter $M = 11$) in simulating realizations of queue lengths, and we calculate the standard deviation for each realization. The STD range in Table 12 shows the minimum and maximum standard deviation for the 10 realizations obtained for each value $\omega$. For each value $\omega$, the STD % Error compares the average STD over the 10 realizations to the STD of 500.5 computed for the observed fillable job openings; and the STD quantile reports the percentage of the 10 realizations for which the STD was less than or equal to the STD of 500.5 computed for observed fillable job openings. Table 12 shows that the STD percent error was lowest and the STD quantile closest to 50% when $\omega = 1.1$. The STD for the observed fillable job openings was not contained in the STD range when $\omega = 1.0$ or $\omega = 4.0$.

The results from Table 12 support the conclusion that $\omega = 1.1$ provides the best fit among the cases considered – with $\omega = 1.2$ not far behind. Additional cases in the range between $\omega = 1.2$ and $\omega = 2.0$ would likely reveal a range of plausible value for $\omega$. Although the best fit appears to occur for a model with a boundary that is graduated to some degree, the comparison in Figure 21 suggests that the estimates of magnitude of censored demands will be similar over the plausible range.

## 8 Bibliography

[1] K. W. Fendick, "Gaussian Fluid Queue with Auutocorrelated Input," *http://arxiv.org/ftp/arxiv/papers/1104/1104.4741.pdf,* 2013.

[2] W. Whitt, Stochastic-Process Limits, Springer, 2002.

[3] J. M. Harrison, Brownian Motion and Stochastic Flow Systems, New York: Wiley, 1985.

[4]  L. Kleinrock, Queueing Systems, Volume 2: Computer Applications, New York: Wiley, 1976.

[5]  G. F. Newell, Applications of Queueing Theory, 2nd Edition, London: Chapman and Hall, 1982.

[6]  R. Cont, S. Stoikov and R. Talreja, "A stochastic model for order book dynamics," *Operations Research,* vol. 58, no. 3, pp. 549-563, 2010.

[7]  W. T. Huh and P. Rusmevichientong, "A nonparametric asymptotic analysis of inventory planning with censored demand," *Mathematics of Operations Research,* vol. 34, no. 1, pp. 103-123, 2009.

[8]  S. A. Conrad, "Sales data and the estimation of demand," *Operations Research Quarterly,* vol. 27, no. 1, pp. 123-127, 1976.

[9]  W. Schnedler, "Likelihood estimation for censored random vectors," *Econometric Reviews,* vol. 24, no. 2, pp. 195-217, 2005.

[10] T. A. Amemiya, Advanced Econometrics, Cambridge, Mass.: Harvard University Press, 1985.

[11] D. R. Cox and P. A. W. Lewis, The Statistical Analysis of Series of Events, London: Methuen, 1966.

[12] A. S. Goldberger, A Course in Econometrics, Cambridge, Mass: Harvard University Press, 1991.

[13] B. Hajek, "A Queue with periodic arrivals and constant service rate," in *Probability, Statistics, and Optimization -- a Tribute to Peter Whittle*, John Wiley and Sons, 1994, pp. 147-158.

[14] J. Kohlas, Stochastic Methods of Operations Research, Cambridge, UK: Cambridge University Press, 1982.

[15] J. F. Hart, E. W. Cheney, C. L. Lawson, H. J. Maehly, C. K. Mesztenyi, J. R. Rice, H. G. Thacher and C. Witzgall, Computer Approximations, New York: Wiley, 1968.

[16] P. Whittle, Prediction and Regulation by Linear Least-Squares Methods, 2nd, Revised ed., Minneapolis: University of Minnesota Press, 1983.

[17] P. Billingsley, Probability and Measure, 2nd Edition, New York: Wiley, 1986.

[18] K. W. Fendick and W. Whitt, "Measurements and approximations to describe the offered traffic and predict the average workload in a single-server queue.," *Proceedings IEEE 77,* pp. 171-194, 1989.

[19] K. W. Fendick, V. R. Saksena and W. Whitt, "Investigating dependence in packet queues with the index of dispersion for work," *IEEE Trans. Commun. 39,* pp. 1231-1243, 1991.

[20] P. B. Stark and R. L. Parker, "Bounded-variable least squares: an algorithm and applications," *Comput Stat., 10, 2,* pp. 129-141, 1995.

[21] J. Mead and R. Renaut, "Least-squares problems with inequality constraints," *Linear Algebra and its Appl., Vol 432, Issue 8, 1,* pp. 1936-1949, 2010.

[22] Wolfstetter and E., Topics in Microeconomics, Cambriidge, U. K.: Cambridge University Press, 2002.

[23] R. W. Wolff, Stochastic Modeling and the Theory of Queues, Upper Saddle River, New Jersey: Prentice Hall, 1986.

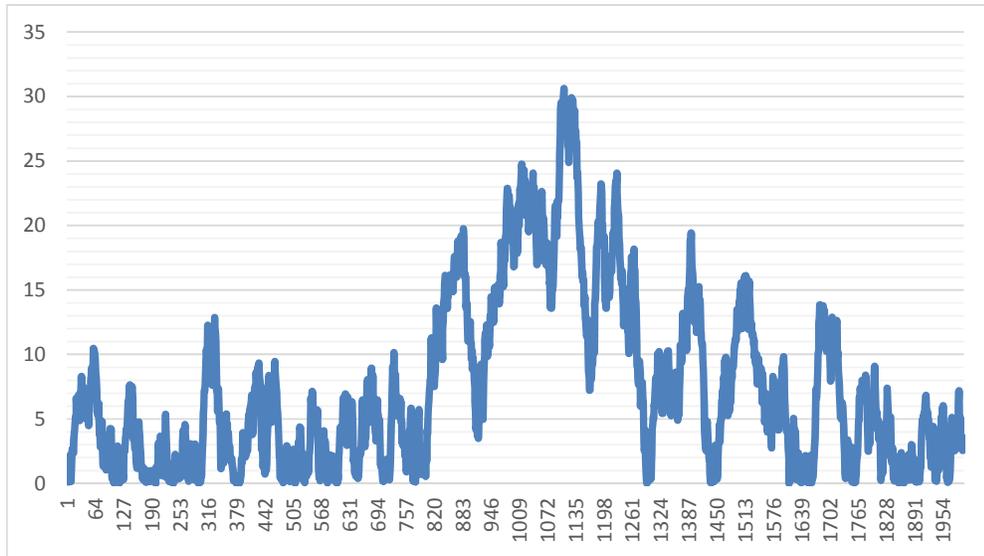

Figure 1. Queue lengths simulated using the M-GMSI model with parameters $M = 3, \rho = -0.1, \theta = 1.1,$ and $\tau = -0.1$.

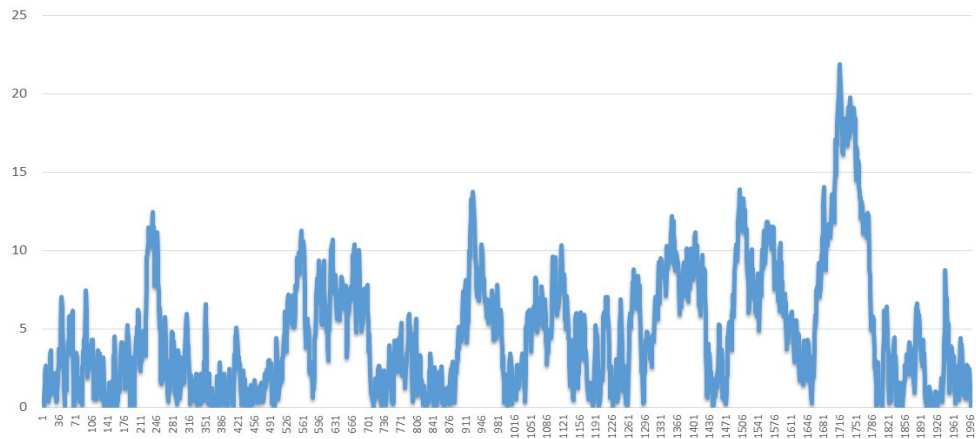

Figure 2. Queue lengths simulated using the M-GMSI model with parameters $M = 3, \rho = -0.1, \theta = 1.1,$ and $\tau = 0.0$.

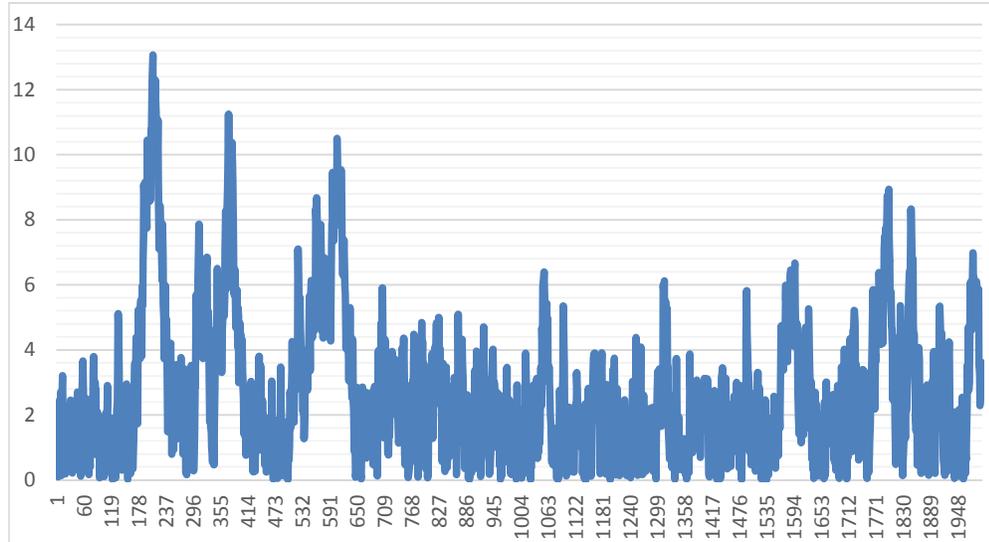

Figure 3. Queue lengths simulated using the M-GMSI model with parameters $M = 3, \rho = -0.1, \theta = 1.1,$ and $\tau = +0.1$.

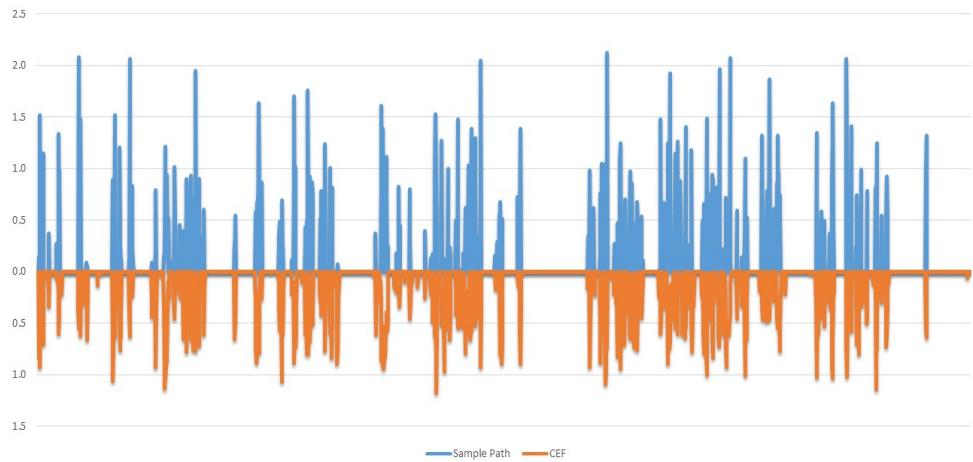

Figure 4. Above the horizontal axis: samples $l_n$ of censored demand from an M-GMSI simulation with parameters $N = 2000, M = 3, \rho = -0.1, \theta = 1.1,$ and $\tau = +0.1$; below the horizontal axis: estimates $\bar{\bar{l}}_n$ of the conditional expectation for censored obtained after 100 steps of DM algorithm using queue lengths from same simulation.

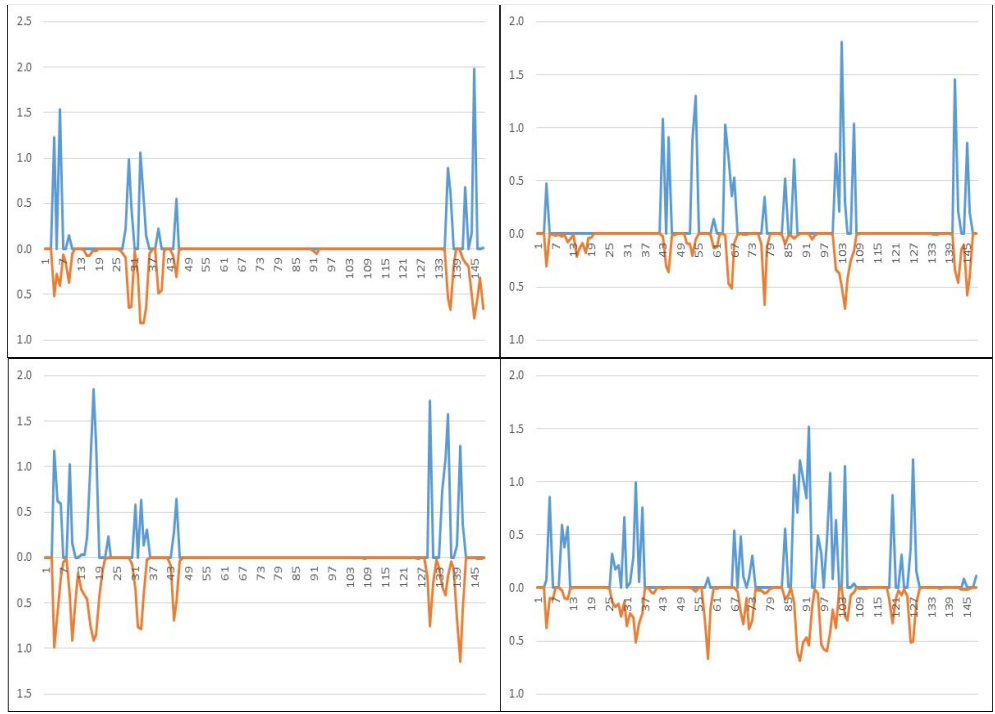

Figure 5. Four realizations from M-GMSI simulations with parameters $N = 2000$, $M = 3, \rho = -0.1, \theta = 1.1,$ and $\tau = +0.1$. Above the horizontal axis for each realization: simulated samples $l_n$ of censored demand; below the horizontal axis for each realization: estimate $\bar{\bar{l}}_n$ of the conditional expectation function for censored demand obtained after 2000 steps of the DM algorithm using queue lengths from same simulation.

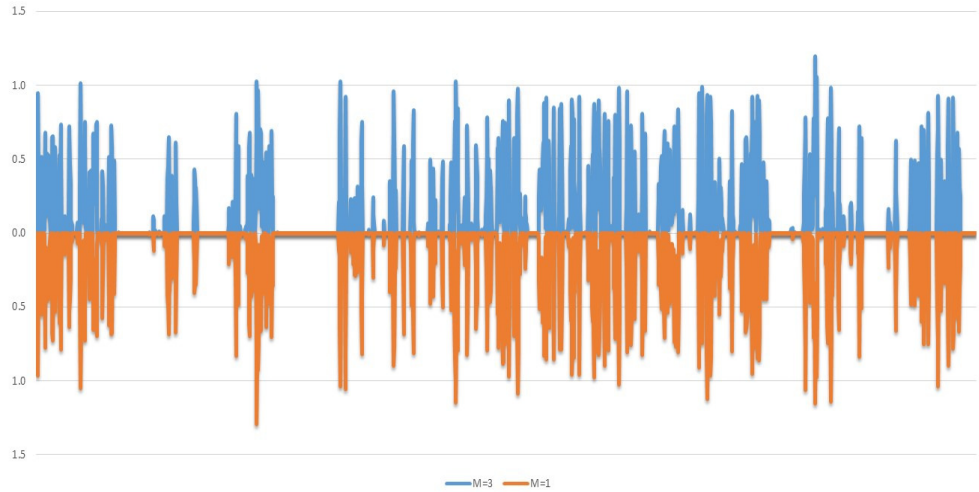

Figure 6. Above the horizontal axis: estimates of conditional expectations $\bar{\bar{l}}_n$ for censored obtained after 100 steps of the DM algorithm assuming the correct value $M = 3$ when sample queue lengths were generated by an M-GMSI simulation with parameters $N = 2000, \ M = 3, \rho = -0.1, \theta = 1.1,$ and $\tau = +0.1$; below the horizontal axis: corresponding estimates $\bar{\bar{l}}_n$ assuming $M = 1$.

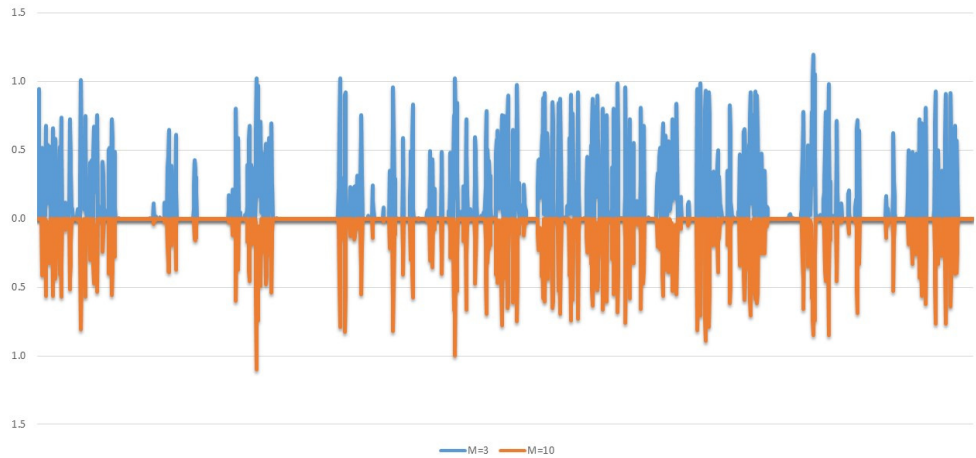

Figure 7. Above the horizontal axis: estimates $\bar{\bar{l}}_n$ of the conditional expectation function for censored obtained after 100 steps of the DM algorithm assuming the correct value $M = 3$ when sample queue lengths were generated by an M-GMSI simulation with parameters $N = 2000, \ M = 3, \rho = -0.1, \theta = 1.1,$ and $\tau = +0.1$; below the horizontal axis: corresponding estimates $\bar{\bar{l}}_n$ assuming $M = 10$.

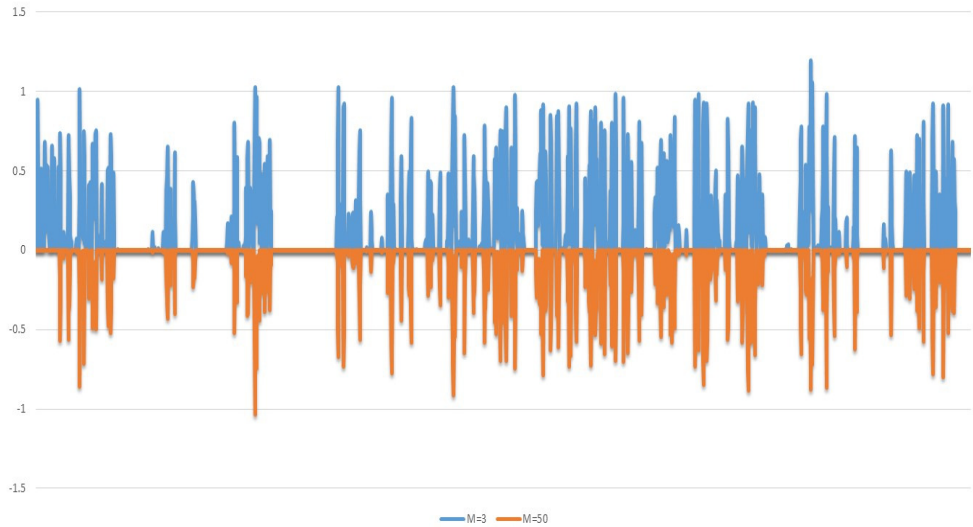

Figure 8. Above the horizontal axis: estimates $\bar{\bar{l}}_n$ of the conditional expectation function for censored obtained after 100 steps of the DM algorithm assuming the correct value $M = 3$ when sample queue lengths were generated by an M-GMSI simulation with parameters $N = 2000$, $M = 3, \rho = -0.1, \theta = 1.1,$ and $\tau = +0.1$; below the horizontal axis: corresponding estimates $\bar{\bar{l}}_n$ assuming $M = 50$.

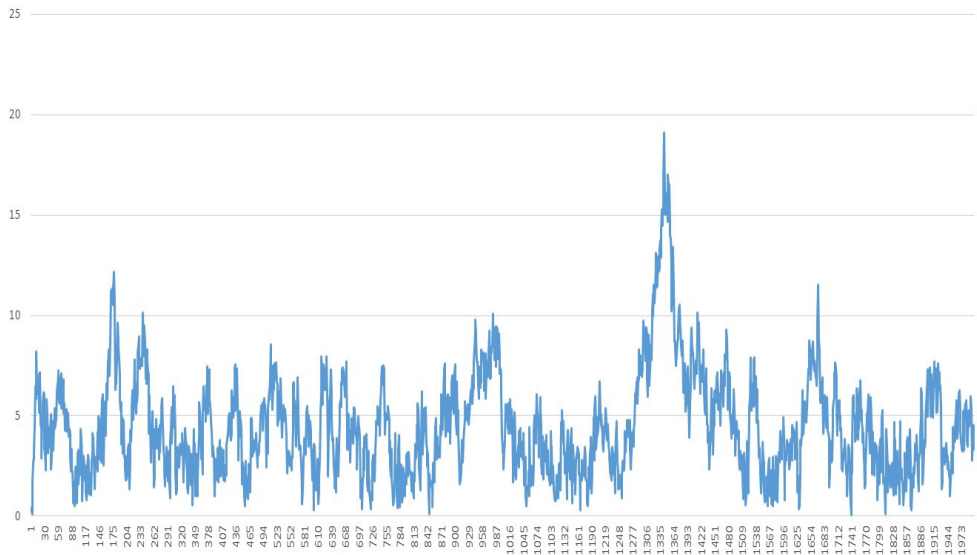

Figure 9. Queue lengths simulated using the $\omega$M-GMSI model with parameters $\omega = 5, M = 3, \rho = -0.1, \theta = 1.1,$ and $\tau = +0.1$.

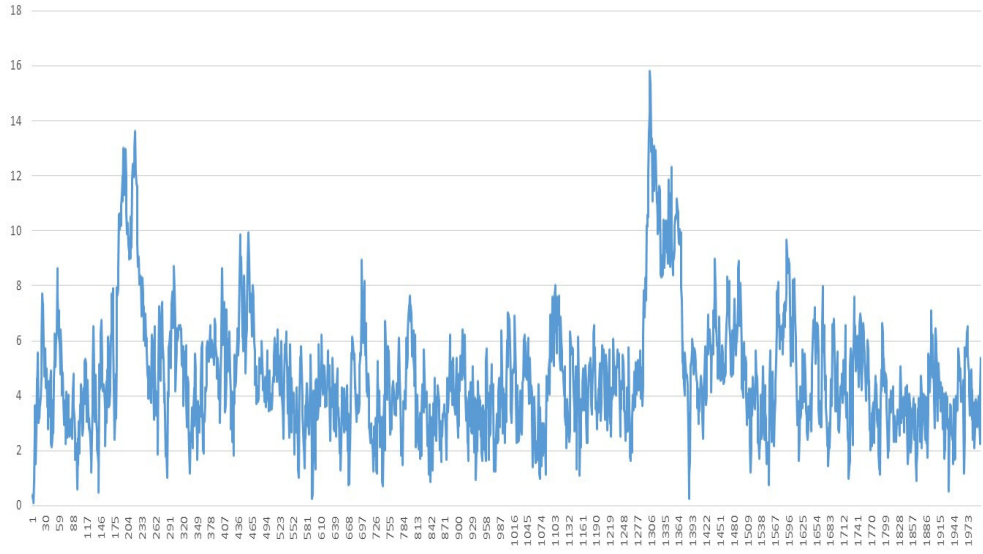

Figure 10. Queue lengths simulated using the $\omega$M-GMSI model with parameters $\omega = 10, M = 3, \rho = -0.1, \theta = 1.1,$ and $\tau = +0.1$.

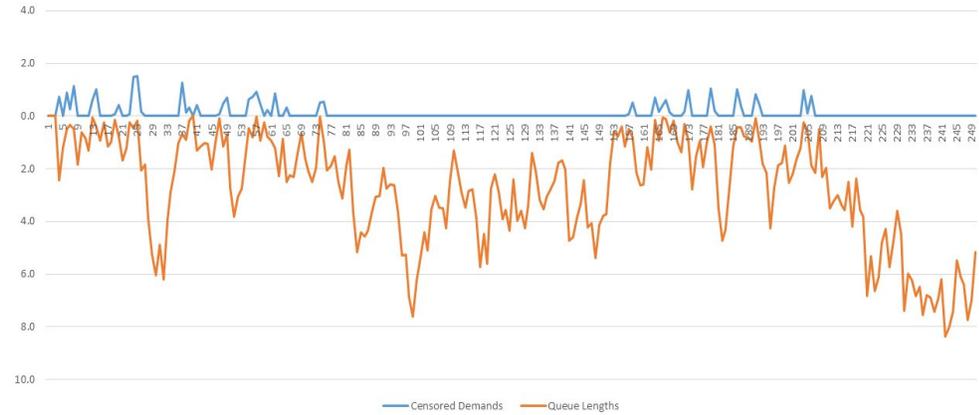

Figure 11. Above the horizontal axis: sample $l_n$ of censored from an $\omega$M-GMSI simulation with parameters $\omega = 1, M = 3, \rho = -0.1, \theta = 1.1,$ and $\tau = +0.1$; below the horizontal axis: queue lengths $q_n$ from same simulation.

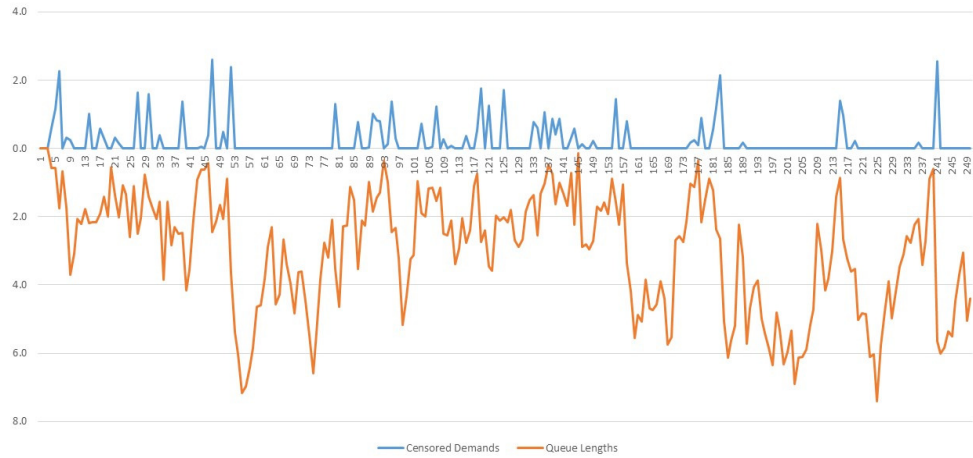

Figure 12.  Above the horizontal axis: sample $l_n$ of censored from an $\omega$M-GMSI simulation with parameters $\omega = 5, M = 3, \rho = -0.1, \theta = 1.1,$ and $\tau = +0.1;$ below the horizontal axis: queue lengths $q_n$ from same simulation.

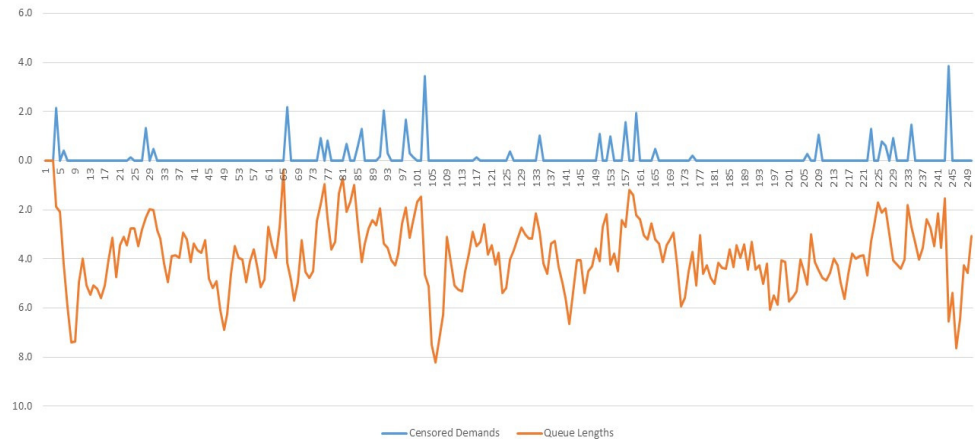

Figure 13.  Above the horizontal axis: sample $l_n$ of censored from an $\omega$M-GMSI simulation with parameters $\omega = 10, M = 3, \rho = -0.1, \theta = 1.1,$ and $\tau = +0.1;$ below the horizontal axis: queue lengths $q_n$ from same simulation.

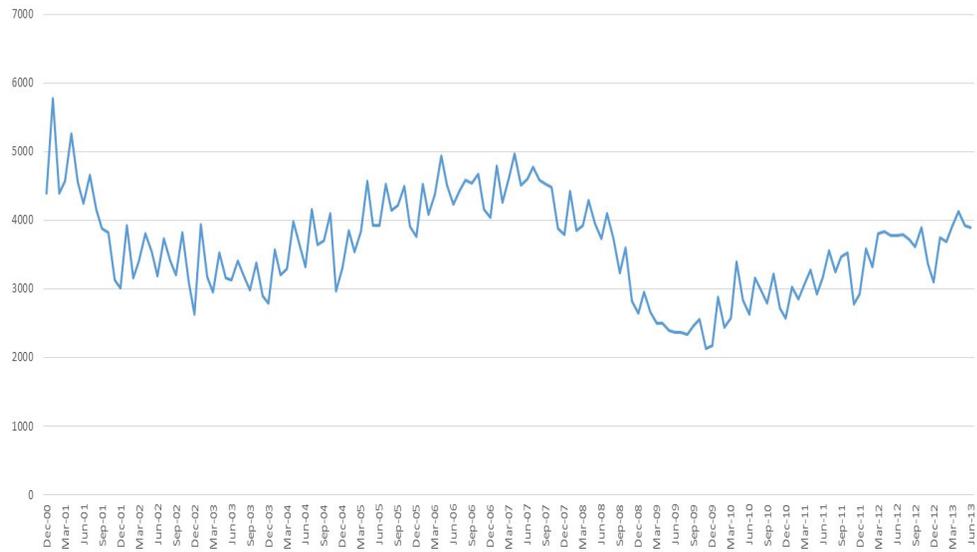

Figure 14. Total non-farm U.S job openings, not seasonally adjusted (source: U.S. Department of Labor: Bureau of Labor Statistics)

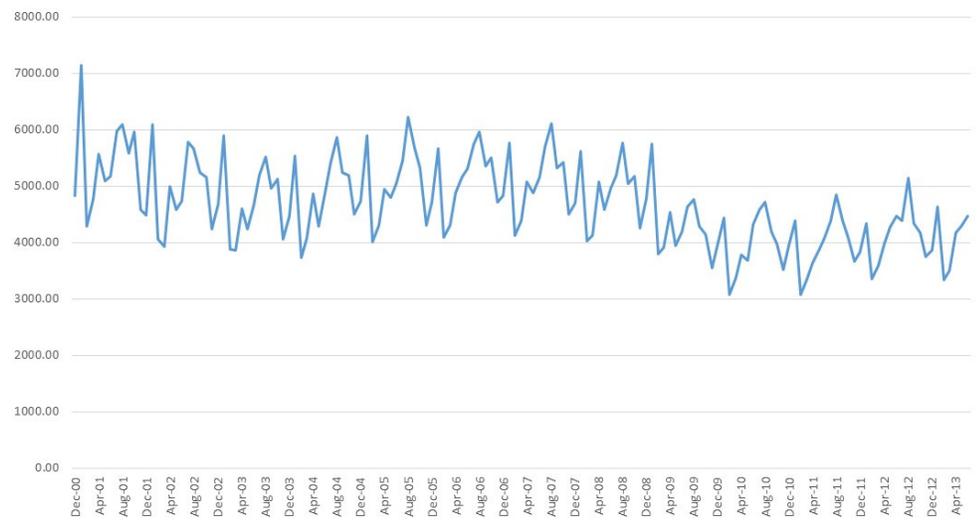

Figure 15. Total non-farm U.S. job separations, not seasonally adjusted (source: U.S. Department of Labor: Bureau of Labor Statistics)

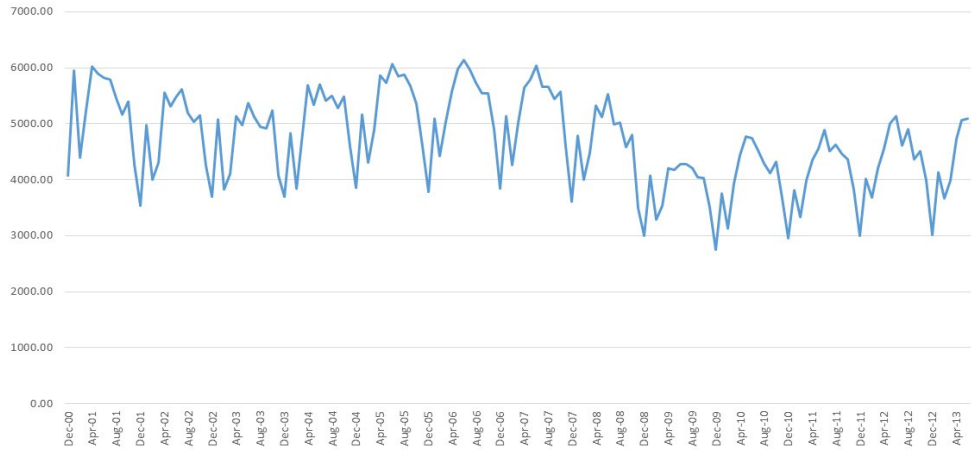

Figure 16. Total non-farm U.S. job hires, not seasonally adjusted (source: U.S. Department of Labor: Bureau of Labor Statistics)

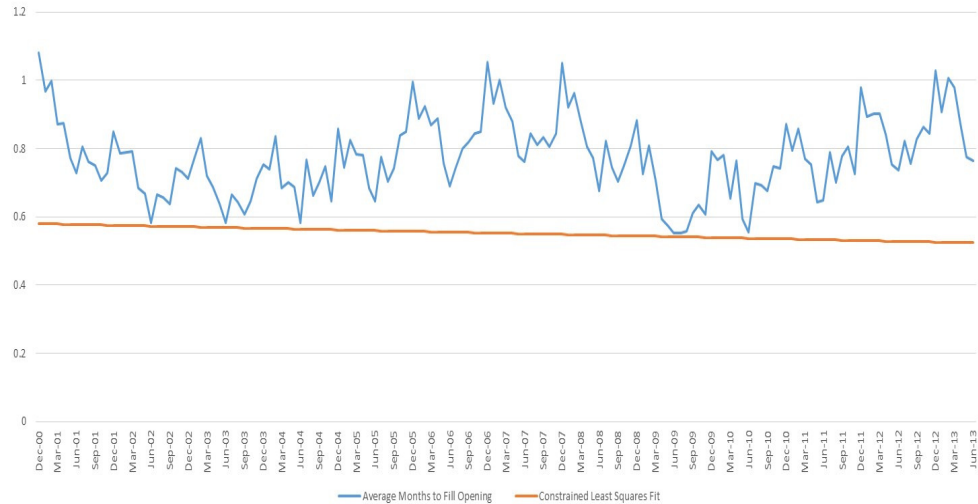

Figure 17. Average months to fill a job opening

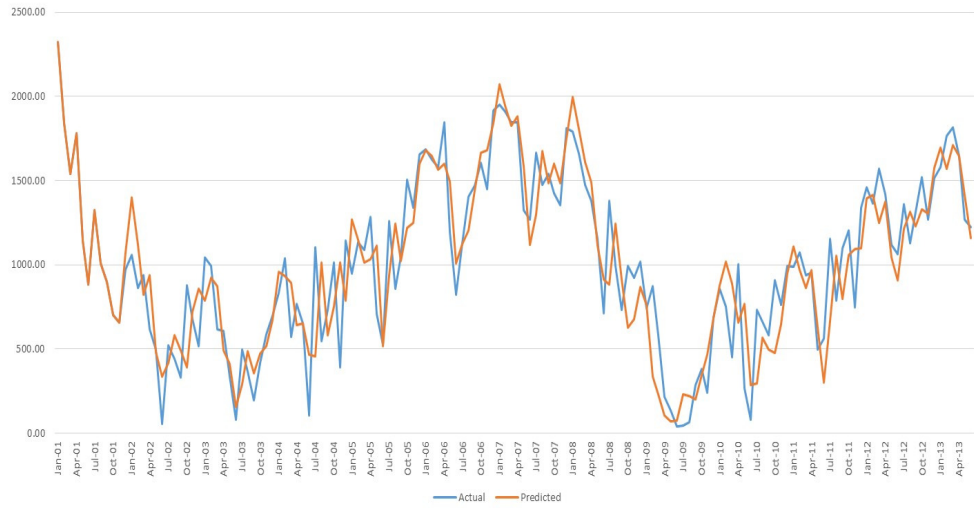

Figure 18. Fillable job openings.

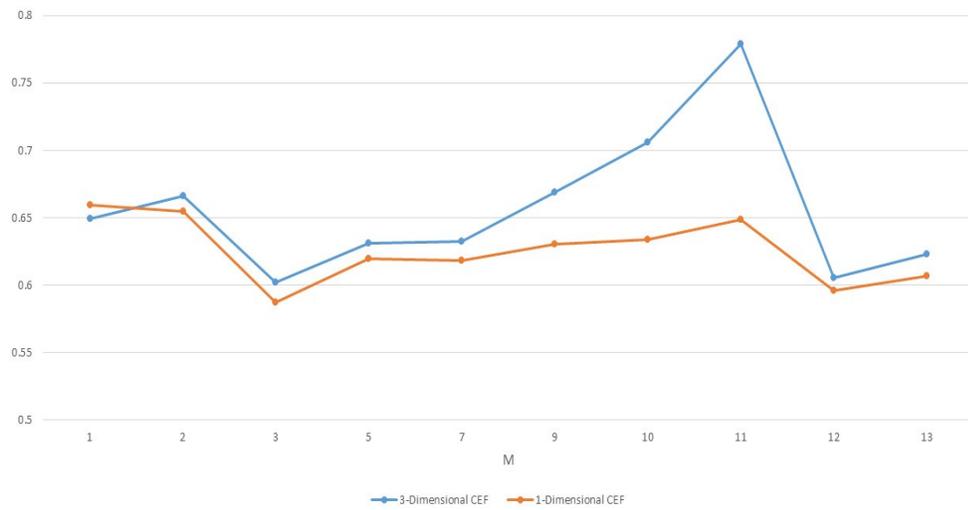

Figure 19. Coefficient of determination ($R^2$) for the one-dimensional and 3-dimensional fits of jobs data

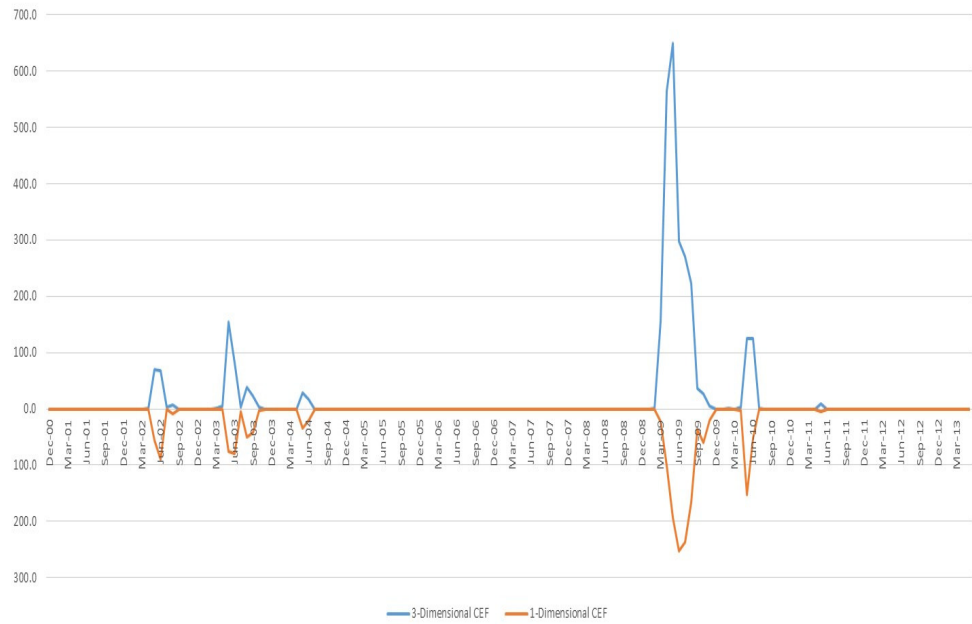

Figure 20. Above the horizontal axis: estimated expected censored demand for jobs using best-fit three-dimensional model; below the horizontal axis: estimated censored demand for jobs using best-fit one-dimensional model.

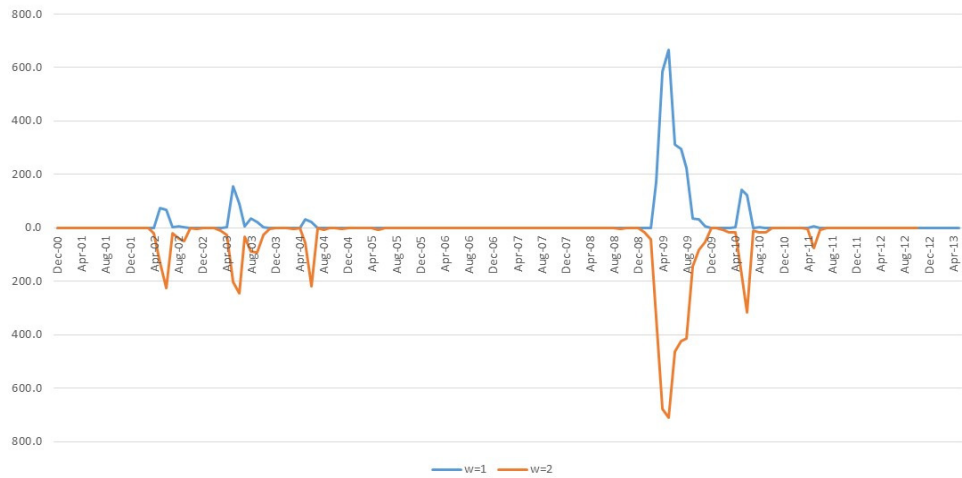

Figure 21. Above the horizontal axis: estimated expected censored demand for jobs using $\omega$M-GMSI modification to three-dimensional model with $M = 11$ and $\omega = 1$; below the horizontal axis: estimated expected censored demand for jobs with $\omega = 2$.

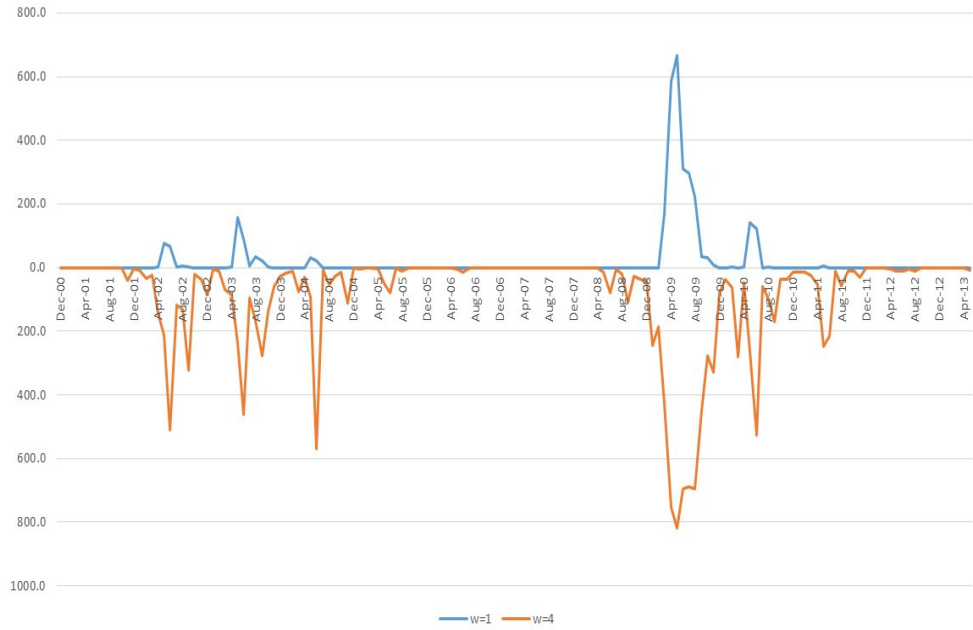

Figure 22. Above the horizontal axis: estimated expected censored demand for jobs using $\omega$M-GMSI modification to three-dimensional model with $M = 11$ and $\omega = 1$; below the horizontal axis: estimated expected censored demand for jobs with $\omega = 4$.

| K | M | ρ | θ | $q_3$ | $z_{1,1}$ | $z_{2,1}$ | $z_{3,1}$ |
|---|---|------|-----|-----|------|------|------|
| 1 | 3 | -0.1 | 1.1 | 0.1 | -0.1 | -0.1 | -0.1 |

Table 1. Common values used for M-GMSI simulations.

| N | $\hat{\rho}$ | $\hat{\theta}$ | $\hat{\tau}$ |
|---|---|---|---|
| 16000 | -0.128 | 1.10 | -0.086 |
| 4000 | -0.094 | 1.10 | -0.078 |
| 149 | -0.139 | 1.14 | -0.049 |
| 149 | -0.091 | 1.26 | -0.170 |
| 149 | -0.228 | 1.07 | -0.166 |
| 149 | -0.266 | 1.02 | -0.215 |
| 149 | -0.109 | 1.19 | -0.068 |
| Weighted Average | -0.123 | 1.10 | -0.086 |

Table 2. Parameter estimates obtained from sample increments of a net-input processes generated by M-GMSI simulations with parameters $M = 3, \rho = -0.1, \theta = 1.1, \ \tau = -0.1$.

| N | $\hat{\rho}$ | $\hat{\theta}$ | $\hat{\tau}$ |
|---|---|---|---|
| 16000 | -0.094 | 1.12 | -0.001 |
| 4000 | -0.103 | 1.03 | -0.001 |
| 149 | +0.025 | 0.94 | -0.012 |
| 149 | -0.102 | 0.89 | -0.043 |
| 149 | -0.061 | 1.01 | +0.010 |
| 149 | -0.028 | 0.87 | +0.034 |
| 149 | -0.261 | 1.03 | -0.025 |
| Weighted Average | -0.095 | 1.10 | -0.001 |

Table 3. Parameter estimates obtained from sample increments of net-input processes generated by M-GMSI simulations with parameters $M = 3, \rho = -0.1, \theta = 1.1,\ and\ \tau = 0.0$.

| N | $\hat{\rho}$ | $\hat{\theta}$ | $\hat{\tau}$ |
|---|---|---|---|
| 16000 | -0.097 | 1.12 | +0.107 |
| 4000 | -0.097 | 1.08 | +0.095 |
| 149 | +0.015 | 1.23 | +0.115 |
| 149 | -0.138 | 0.94 | +0.087 |
| 149 | -0.075 | 0.93 | +0.021 |
| 149 | -0.087 | 1.00 | +0.108 |
| 149 | +0.004 | 0.98 | +0.078 |
| Weighted Average | -0.093 | 1.11 | +0.103 |

Table 4. Parameter estimates obtained from sample increments of net-input processes generated by an M-GMSI simulations with parameters $M = 3, \rho = -0.1, \theta = 1.1,$ and $\tau = +0.1$

| M | $\hat{\rho}$ | $\hat{\theta}$ | $\hat{\tau}$ |
|---|---|---|---|
| 1 | -0.097 | 1.11 | +0.101 |
| 2 | -0.097 | 1.11 | +0.102 |
| 3 | -0.097 | 1.12 | +0.107 |
| 4 | -0.097 | 0.97 | +0.070 |
| 10 | -0.097 | 0.82 | +0.025 |
| 50 | -0.097 | 0.79 | +0.005 |

Table 5. Parameter estimates obtained from sample increments of a net-input process generated by an M-GMSI simulation with parameters $N = 16000$, $M = 3, \rho = -0.1, \theta = 1.1,$ and $\tau = +0.1$.

| Steps | $\hat{\rho}$ | $\hat{\theta}$ | $\hat{\tau}$ |
|---|---|---|---|
| 0 | +0.001 | 1.02 | +0.116 |
| 1 | -0.083 | 1.09 | +0.109 |
| 2 | -0.098 | 1.10 | +0.106 |
| 3 | -0.108 | 1.10 | +0.098 |
| 4 | -0.109 | 1.12 | +0.102 |
| 5 | -0.113 | 1.13 | +0.102 |
| 6 | -0.111 | 1.12 | +0.099 |
| Average (steps 1-6) | -0.108 | 1.11 | +0.101 |

Table 6. Parameter estimates after successive steps of the DM algorithm from sample queue lengths generated by an M-GMSI simulation with parameters $N = 8000$, $M = 3, \rho = -0.1, \theta = 1.1,$ and $\tau = +0.1$.

| Steps | N | $\tau$ | $\overline{\hat{\rho}}$ | $\overline{\hat{\theta}}$ | $\overline{\hat{\tau}}$ |
|---|---|---|---|---|---|
| 5 | 8000 | -0.1 | -0.108 | 1.12 | -0.091 |
| 5 | 8000 | -0.1 | -0.120 | 1.10 | -0.100 |
| 5 | 8000 | -0.1 | -0.100 | 1.11 | -0.090 |
| 5 | 8000 | -0.1 | -0.120 | 1.13 | -0.098 |
| 100 | 149 | -0.1 | -0.112 | 1.29 | -0.076 |
| 5 | 8000 | 0.0 | -0.116 | 1.14 | -0.007 |
| 5 | 8000 | +0.1 | -0.104 | 1.14 | +0.103 |
| 5 | 8000 | +0.1 | -0.107 | 1.11 | +0.101 |
| 5 | 8000 | +0.1 | -0.101 | 1.09 | +0.099 |
| 5 | 8000 | +0.1 | -0.11 | 1.04 | +0.093 |
| 100 | 149 | +0.1 | -0.138 | 0.95 | +0.089 |

Table 7. Parameter estimates after the given number of steps of the DM algorithm from sample queue lengths generated by M-GMSI simulations with parameters $M = 3, \rho = -0.1,$ and $\theta = 1.1$.

| M | $\overline{\hat{\rho}}$ | $\overline{\hat{\theta}}$ | $\overline{\hat{\tau}}$ | $R^2$ |
|---|---|---|---|---|
| 1 | -0.116 | 1.10 | +0.086 | 0.88 |
| 2 | -0.115 | 1.11 | +0.094 | 0.88 |
| 3 | -0.112 | 1.16 | +0.098 | 0.88 |
| 4 | -0.092 | 0.94 | +0.065 | 0.88 |
| 10 | -0.075 | 0.78 | +0.026 | 0.88 |
| 50 | -0.065 | 0.67 | +0.006 | 0.88 |

Table 8. Parameter estimates after 50 steps of the DM algorithm (assuming different values $M$) from sample queue lengths generated by an M-GMSI simulation with parameters $N = 8000, \ M = 3, \rho = -0.1, \theta = 1.1,$ and $\tau = +0.1$.

| $\omega$ | $\overline{\hat{\rho}}$ | $\overline{\hat{\theta}}$ | $\overline{\hat{\tau}}$ | $R^2$ | STD % Error |
|---|---|---|---|---|---|
| 1 | -0.009 | 1.12 | +0.125 | 0.85 | 209.2% |
| 3 | -0.048 | 1.11 | +0.115 | 0.85 | 33.7% |
| 5 | -0.094 | 1.07 | +0.097 | 0.85 | 7.8% |
| 7 | -0.141 | 1.04 | +0.078 | 0.85 | 8.4% |

Table 9. Parameter estimates after 50 steps of the ωDM algorithm (assuming different values $\omega$) from sample queue lengths generated by an ωM-GMSI simulation with parameters $N = 8000, M = 3, \omega = 5.0, \rho = -0.1, \theta = 1.1,$ and $\tau = +0.1$.

| M | $\bar{\bar{\rho}}$ | $\bar{\bar{\theta}}$ | $\bar{\bar{\tau}}$ | $R^2$ |
|---|---|---|---|---|
| 1 | -0.030 | 0.202 | 0.049 | 0.66 |
| 2 | -0.022 | 0.166 | 0.032 | 0.65 |
| 3 | -0.018 | 0.084 | 0.004 | 0.59 |
| 5 | -0.016 | 0.137 | 0.012 | 0.62 |
| 7 | -0.015 | 0.117 | 0.007 | 0.62 |
| 9 | -0.012 | 0.116 | 0.007 | 0.63 |
| 11 | -0.013 | 0.136 | 0.008 | 0.65 |
| 13 | -0.016 | 0.075 | 0.002 | 0.61 |

Table 10. Parameter estimates based on a 1-dimensional fit ($K = 1$) of jobs data by the DM algorithm after 500 steps for varying values for $M$.

| M | $\bar{\bar{\rho}}$ | $\bar{\bar{\theta}}$ | $\bar{\bar{\tau}}$ | $R^2$ |
|---|---|---|---|---|
| 1 | -0.030 | 0.202 | 0.049 | 0.65 |
| 2 | -0.023 | 0.166 | 0.033 | 0.66 |
| 3 | -0.027 | 0.084 | 0.003 | 0.60 |
| 5 | -0.019 | 0.138 | 0.012 | 0.63 |
| 7 | -0.020 | 0.120 | 0.007 | 0.63 |
| 9 | -0.018 | 0.120 | 0.007 | 0.67 |
| 11 | -0.029 | 0.162 | 0.006 | 0.78 |
| 13 | -0.014 | 0.079 | 0.002 | 0.62 |

Table 11. Parameter estimates based on a 3-dimensional fit ($K = 3$) of jobs data by the DM algorithm after 500 steps with varying values for $M$.

| | | | | | STD | | |
|---|---|---|---|---|---|---|---|
| $\omega$ | $\bar{\bar{\rho}}$ | $\bar{\bar{\theta}}$ | $\bar{\bar{\tau}}$ | $R^2$ | Range | % Error | Quantile |
| 1.0 | -0.029 | 0.16 | 0.007 | 0.779 | 567.5-909.9 | 30% | 0% |
| 1.1 | -0.031 | 0.16 | 0.006 | 0.777 | 453.0-830.9 | 14% | 20% |
| 1.2 | -0.034 | 0.16 | 0.006 | 0.772 | 473.0-746.0 | 21% | 10% |
| 2.0 | -0.052 | 0.17 | 0.004 | 0.766 | 492.2-1058.8 | 33% | 10% |
| 4.0 | -0.110 | 0.19 | -0.006 | 0.731 | 595.2-997.8 | 53% | 0% |

Table 12. Estimates of parameters and $R^2$ based on a 3-dimensional fit ($K = 3$) of jobs data by the ωDM algorithm with $M = 11$ and varying values $\omega$. Estimates of standard deviation of queue lengths ("STD") for each value $\omega$ based on 10 independent ωM-GMSI realizations. The "% Error" and "Quantile" for STD are relative to the standard deviation of observed queue lengths $\{q_n\}$ of 500.5.